\documentclass{article}[11pt]

\usepackage{amssymb}
\usepackage{amsmath}
\usepackage{amscd}
\usepackage{amsthm}
\usepackage{color}
\usepackage{graphicx}
\usepackage{subfig}
\usepackage{url}
\usepackage{mdwlist}
\usepackage{tabls}
\usepackage{array}
\usepackage{algorithmic}
\usepackage{algorithm}

\usepackage{pgf}
\usepackage{tikz}

\usepackage{verbatim} 

\newcommand*{\CopyCounter}[2]{%
  \expandafter\def\csname c@#2\endcsname{\csname c@#1\endcsname}%
  \expandafter\def\csname p@#2\endcsname{\csname p@#1\endcsname}%
  \expandafter\def\csname the#2\endcsname{\csname the#1\endcsname}}

\numberwithin{Theorem}{section}
\CopyCounter{Theorem}{Proposition}
\CopyCounter{Theorem}{ProposedProblem}
\CopyCounter{Theorem}{Property}
\CopyCounter{Theorem}{Claim}
\CopyCounter{Theorem}{Lemma}
\CopyCounter{Theorem}{Corollary}
\CopyCounter{Theorem}{Conjecture}
\CopyCounter{Theorem}{Definition}
\CopyCounter{Theorem}{Example}
\CopyCounter{Theorem}{Remark}
\CopyCounter{Theorem}{Question}
\CopyCounter{Theorem}{Condition}
\CopyCounter{Theorem}{Criterion}
\CopyCounter{Theorem}{Observation}
\theoremstyle{plain}

\theoremstyle{definition}

\newtheorem{remark}[Remark]{Remark}

\newcommand{\bm}[1]{\mbox{\boldmath${#1}$}}

\newcommand{\domain}{\Omega}
\newcommand{\cdomain}{\widebar{\Omega}}
\newcommand{\boundary}{\partial \domain}

\newcommand{\xbar}{\bm{\bar{x}}}

\newcommand{\x}{\bm{x}}

\newcommand{\ba}{\bm{a}}

\newcommand{\y}{\bm{y}}

\newcommand{\widebar}{\overline}
\newcommand{\R}{\bm{R}}

\newcommand{\J}{{\cal{J}}}

\newcommand{\fB}{\bm{f}}

\newcommand{\drop}[1]{}

\newif\ifnever\neverfalse

\makeatletter
\newcommand{\rmnum}[1]{\romannumeral #1}
\newcommand{\Rmnum}[1]{\expandafter\@slowromancap\romannumeral #1@}
\makeatother

\newcommand{\marginfix}{
\setlength{\parskip}{0.01cm}
\setlength{\textwidth}{6.0in}
\setlength{\oddsidemargin}{-0.0 in}
\setlength{\evensidemargin}{0.0 in}
\setlength{\topmargin}{-0.5in}
\setlength{\textheight}{9.0 in}
}

\marginfix

  \renewenvironment{thebibliography}[1]{%
    \begin{oldthebibliography}{#1}%
      \setlength{\parskip}{.3ex}%
      \setlength{\itemsep}{.3ex}%
  }%
  {%
    \end{oldthebibliography}%
  }

\let\oldmarginpar\marginpar
\renewcommand\marginpar[1]{\-\oldmarginpar[\raggedleft\footnotesize #1]{\raggedright\footnotesize #1}}

\newcommand{\todo}[1]{{\large $\spadesuit$ {\bf #1} $\spadesuit$ }}

\begin{document}


{\Large{\bf
\centerline{A fast implicit method for time-dependent Hamilton-Jacobi PDEs.}
}}

\vspace*{.2in}
{\Large
\centerline{Alexander Vladimirsky\footnote{
\mbox{Department of Mathematics, Cornell University, Ithaca, NY 14853}
}
 and Changxi Zheng\footnote{
Department of Computer Science, Columbia University, New York, NY 10027}
}

}

\vspace*{.1in}
\begin{abstract}
\noindent
We present a new efficient computational approach for time-dependent
first-order Hamilton-Jacobi-Bellman PDEs.  Since our method is based on
a time-implicit Eulerian discretization, the numerical scheme is
unconditionally stable, but discretized equations for each
time-slice are coupled and non-linear.
We show that the same system can be re-interpreted as
a discretization of a static Hamilton-Jacobi-Bellman PDE on
the same physical domain.  The latter was shown
to be ``causal'' in \cite{VladTimeD}, making fast (non-iterative)
methods applicable.
The implicit discretization results in higher computational cost per time
slice compared to the explicit time marching.
However, the latter is subject to a CFL-stability condition,
and the implicit approach becomes significantly more efficient whenever
the accuracy demands on the time-step are less restrictive than
the stability.  We also present a hybrid
method, which aims to combine the advantages of both the explicit and
implicit discretizations.
We demonstrate the efficiency of our approach using
several examples in optimal control of isotropic fixed-horizon
processes.
\end{abstract}

\vspace*{.1in}

\vspace*{.1in}

\section{Introduction.}
\label{s:intro}
For evolutive partial differential equations, explicit time marching
provides a popular and conceptually simple computational approach.
However, its main drawback is due to the CFL-type stability conditions,
which often result in a choice of time-steps significantly smaller
than what could be expected based on the accuracy requirements.
On the other hand, implicit time schemes are usually
unconditionally stable, but for non-linear PDEs they result in
a discretized system of coupled nonlinear equations, which need to be
solved at each time step.  This task can be prohibitively expensive,
if the system is solved using iterative methods.  A similar challenge
also exists for discretizations of static nonlinear boundary value problems.
Nevertheless, fast (non-iterative) methods have been developed
for a class of such static PDEs (first-order Hamilton-Jacobi-Bellman equations).

In this paper we develop efficient implicit (and ``hybrid'' explicit-implicit)
methods for {\em time-dependent} HJB equations by re-interpreting each time-slice as a
discretization of an auxiliary {\em static} HJB problem.

We start by comparing explicit, implicit, and hybrid approaches in a simpler setting
of linear 1D advection equations with  non-homogeneous advection speeds (section \ref{s:advection}).
We provide a quick review of time-dependent and static HJB PDEs
in the context of deterministic optimal control theory in section \ref{s:HJB_intro}.
We then describe the standard first-order accurate
explicit and implicit discretizations for the evolutive case
and explain the connection to the discretizations
of static HJB equations (section \ref{s:discretize}).
Non-iterative methods for a relevant class of the latter are reviewed
in section \ref{ss:labelset} and our algorithm for the time-dependent
case is summarized in section \ref{ss:algorithms}.

It is worth noting that the convergence and unconditional stability of implicit schemes
for Hamilton-Jacobi PDEs are well known.  Indeed, the implicit discretization used in this paper
falls into the class of schemes analyzed by Souganidis \cite{Souganidis85} back in 1985.
But such schemes were previously considered impractical -- and the challenge of solving
a coupled non-linear system of discretized equation in each time-slice was only a part of
the explanation for this.
Implicit methods were also considered unnecessary
because of the perception that, for first-order PDEs, the accuracy requirements
on time-steps are typically at least as restrictive as the CFL stability conditions.
Such arguments are often made even for hyperbolic conservation
laws; e.g., in \cite{Kroner}, the analysis of 1D advection equation with
constant coefficients 
is used to show that the errors resulting from
implicit methods are larger than those for explicit methods,
when the same time-step is used in both.
However, 
in realistic problems
characteristic speeds might vary significantly
throughout the computational domain.
In recognition of this, several implicit and hybrid methods were 
developed for
linear and quasi-linear problems with stiff boundary layers \cite{CollinsColellaGlaz, ORourkeSahota}.
Our numerical results in section \ref{s:experiments} demonstrate
that implicit and hybrid (implicit/explicit) methods
are similarly significantly better for a class of (fully-nonlinear)
Eikonal problems with strong time-space inhomogeneities.
All considered examples can be naturally interpreted as dynamic programming equations
for the value functions of fixed-horizon isotropic optimal control problems.

We conclude by discussing 
possible extensions
and directions for future work in section \ref{s:conclusions}.

\section{Explicit, implicit, and hybrid methods for advection in 1D.}
\label{s:advection}
Implicit schemes for non-linear evolutive hyperbolic equations have been long considered non-competitive primarily for two reasons: (\rmnum{1}) implicit time schemes for
non-linear PDEs result in a discretized system of coupled nonlinear equations,
which is generally expensive to solve, (\rmnum{2}) for first-order PDEs, the
accuracy requirement on time-steps are typically as restrictive as CFL
stability conditions, and therefore the advantage of using large time-steps in
implicit schemes becomes insignificant. The latter argument is frequently made
even for advection equations \cite{Kroner}. 
Postponing the discussion of argument (\rmnum{1}) until section \ref{s:discretize},
here we focus on the argument (\rmnum{2}) for simpler linear hyperbolic PDEs in 1D.
We show that, for a wide range of strongly non-homogeneous advection speeds,
the ``computational cost per accuracy'' of implicit schemes is lower than that of explicit schemes.
This analysis remains largely valid even for the more general PDEs considered in the rest of this paper, but with a few caveats due to the higher dimensional state space and the non-linearity of HJB equations.
Some of the material in this section is fairly standard, but we cover it from a different perspective, to set up the context for sections \ref{s:HJB_intro} and \ref{s:discretize}.

In 1D the direct relation between HJB PDEs and hyperbolic conservation laws (HCLs)
is well-known: if $v(x,t)$ is a $C^2$ solution to a HJB PDE $v_t + H(v_x,x,t)=0$, then $w(x,t)=v_x(x,t)$ solves an HCL
$w_t + \left[ H(w,x,t) \right]_x =0$.
At the same time, the non-divergence-form advection problem
\begin{eqnarray}
\label{eq:advect}
v_t + f(x,t) v_x = 0, && \qquad \forall x \in \domain = (0,1), \, t \geq 0;\\
\nonumber
v(x,0) = \alpha(x); \quad
v(0,t) = \beta(t); \quad
f(x,t) > 0; &&
\end{eqnarray}
can be also re-interpreted as a HJB equation corresponding to an ``optimal'' control problem with no running cost and a singleton set of control values; see also Remark \ref{rem:HJB_as_advect}.

The spatial derivatives can be approximated using appropriate divided differences, resulting in
a semi-discretization of~\eqref{eq:advect}.
Using 1D grid functions $V_i(t) \approx v(hi,t)$ and $F_i=f(hi,t)$ with the spatial gridsize $h=1/M$,
we approximate this PDE with a system of ODEs
\begin{eqnarray}\label{eq:semi-discr}
    \frac{d}{dt}V_i &=& -F_i(t)\frac{V_i-V_{i-1}}{h},\\
    \nonumber
    V_i(0) &=& \alpha(hi), \qquad \qquad \qquad \qquad i=1\dots M;\\
    \nonumber
    V_0(t) &=& \beta(t),
\end{eqnarray}
where $v_x$ is approximated using the first-order upwind scheme.
The usual fully discrete explicit and implicit schemes are respectively the forward and
backward Euler time discretizations of~\eqref{eq:semi-discr}. In particular,
defining $\lambda^{n}_i = \frac{k}{h}f(hi,nk)$ and using the time-step $k>0$,
the explicit scheme computes $V^{n+1}_i \approx v(hi,(n+1)k)$ using
\begin{equation}\label{eq:exp_scheme_1d}
    V^{n+1}_i = \lambda^{n}_i V^n_{i-1} + (1-\lambda^{n}_i) V^n_i,
\end{equation}
and the implicit scheme computes $V^{n+1}_i$ using
\begin{equation}\label{eq:imp_scheme_1d}
    V^{n+1}_i = \frac{1}{\lambda^{n+1}_i+1} V^n_{i} + \frac{\lambda^{n+1}_i}{\lambda^{n+1}_i+1} V^{n+1}_{i-1}.
\end{equation}
The time step of the explicit scheme is restricted with a CFL stability condition
$\max_{i,n} \lambda_i^n \leq 1$ or, equivalently, $k \leq \hat{k} = h / \hat{f},$ where $ \hat{f} = \max_{x,t} f(x, t).$
So, the total cost of computing the solution up to the time $T = Nk$ is $O(M N) = O(T M^2 \hat{f}).$
In contrast, the implicit scheme is unconditionally stable, resulting in a potentially much smaller computational cost of $O(MT/k)$ when $k$ is significantly larger than $\hat{k}$.

\begin{remark}\label{rem:lin_imp_cost}
The above cost comparison is based on simply counting the number of gridpoint values that need to be computed up to the time $T$.  We note that the linearity of the PDE  and the simple 1D flow nature of the problem ensure that in the implicit scheme $V^{n+1}_i$ can be computed directly/sequentially based on the already-known grid value $V^{n+1}_{i-1}$.  Thus, the cost of each gridpoint update is the same in both explicit and implicit schemes
(equations~\eqref{eq:exp_scheme_1d} and \eqref{eq:imp_scheme_1d}).  This is quite different from the general non-linear case in $R^d$, where implicit updates are typically more expensive and an even larger $k$ is needed to realize any computational savings.  See also Remark \ref{rem:cost_update}.
\end{remark}

\begin{figure}[ht]
    \centering
    \subfloat[Explicit scheme is better.]{\includegraphics[width=0.425\columnwidth]{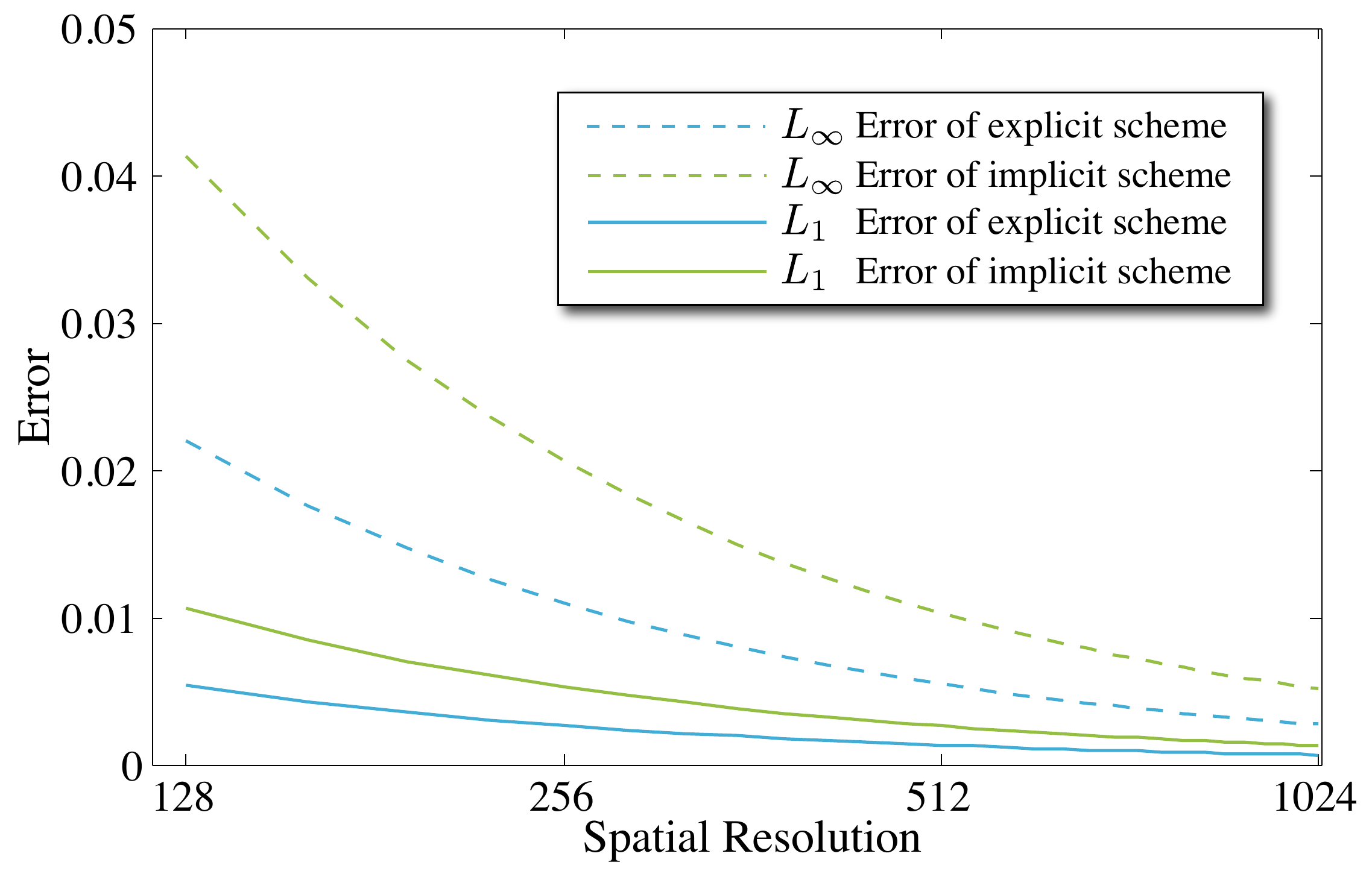}\label{fig:exp_wins}}
    \subfloat[Implicit scheme is better.]{\includegraphics[width=0.41\columnwidth]{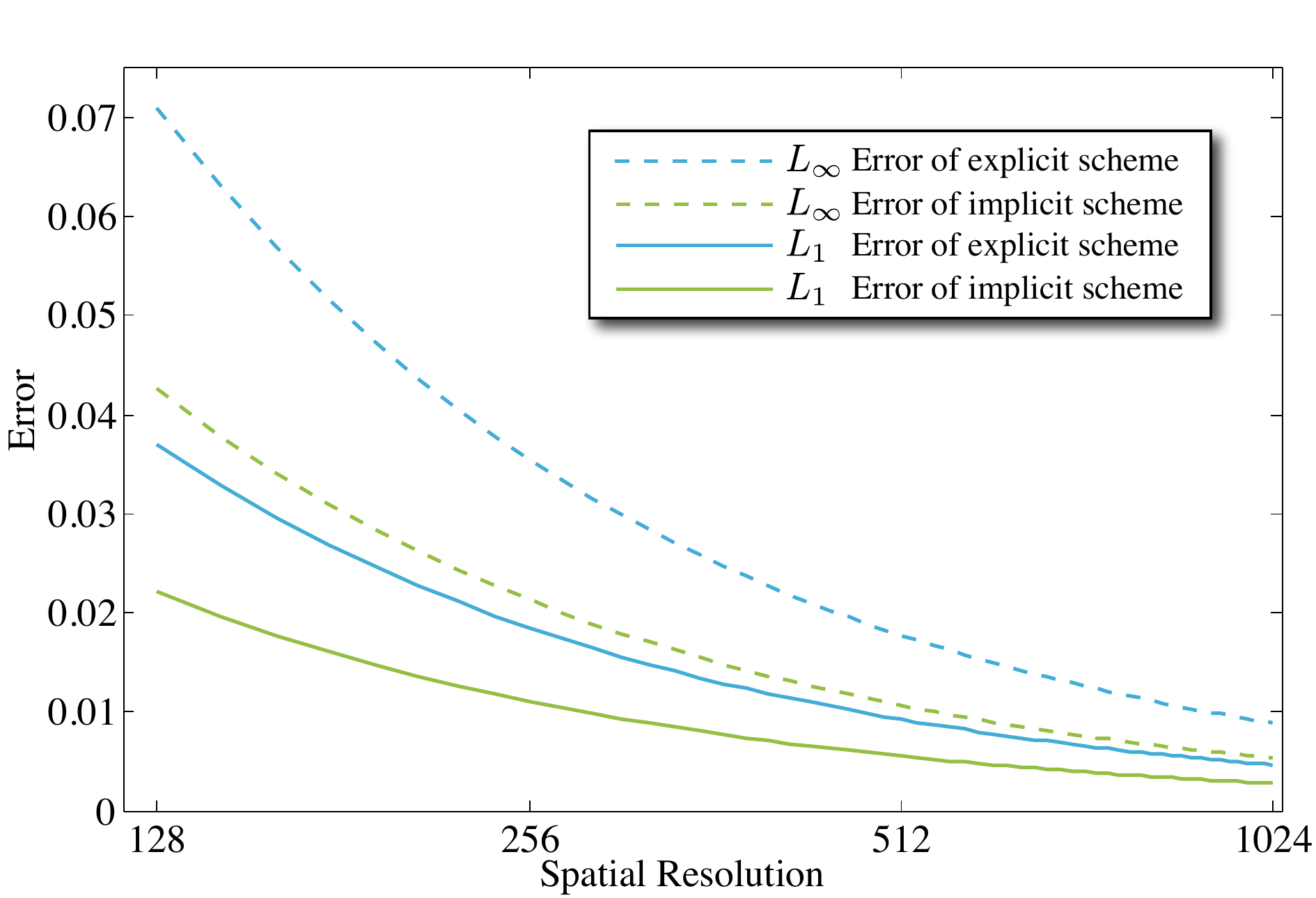}\label{fig:imp_wins}}\;\;
    \subfloat[Stencils]{\includegraphics[width=0.14\columnwidth]{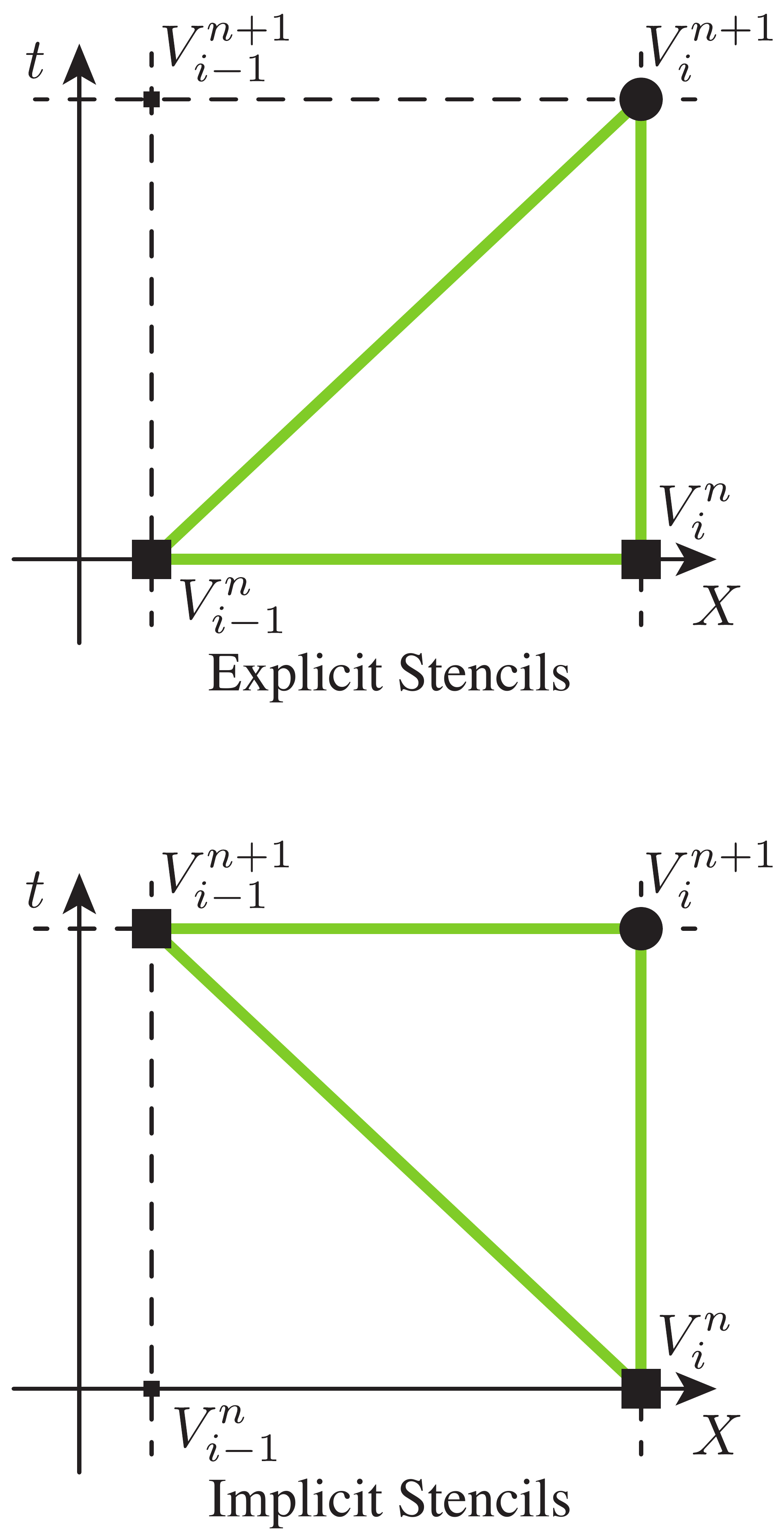}}
    \caption{\textbf{Accuracy comparison of explicit and implicit schemes:}
    $L_1$ and $L_\infty$ errors in the approximate solution at $T=1.5$.
    Both explicit and implicit schemes use the same time-step
    determined by the CFL condition at the given resolution. 
    Left to right:
    (a) $v_t + \frac{1}{2x+1}v_x = 0$ with $v(0,t)=e^{-t}$, and the analytic solution $v(x,t) = e^{x^2+x-t}$;
    (b) $v_t + v_x = 1.5e^{x+t}$ with $v(0,t) = e^{t}$, and the analytic solution $v(x,t) = e^{x+t}$;
    (c) different stencils used in both schemes.
    \label{fig:imp_exp_win}}
\end{figure}

For a fair comparison, any computational savings should be also balanced against the method's accuracy.
E.g., if $f$ is constant then the time step $\hat{k} = h/f$ is obviously optimal for the explicit scheme, resulting in the exact solution since $v(x_i, t_{n+1}) = v(x_i - k f, t_{n+1}-k) = v(x_{i-1}, t_n)$.  For this simple (constant advection speed) example, it is well-known that\\
$\bullet \,$ the accuracy of the explicit scheme actually deteriorates if we insist on using some smaller time-step $k < \hat{k}$;\\
$\bullet \,$  for any $k \leq \hat{k}$, the errors in the implicit scheme will be larger than in the explicit scheme.\\
This is best understood by analyzing the {\em modified equations} for each scheme, which are PDEs with $h$-dependent coefficients, satisfied by the numerical solution to a higher-order of accuracy than the original equation \eqref{eq:advect}.
For the spatially non-homogeneous case $f(x,t) = f(x)$, the Taylor series expansion shows that the modified equations are respectively
\begin{eqnarray}
    v_t + f v_x &=& \frac{1}{2} h f v_{xx}(1-\lambda) - \frac{1}{2}h \lambda f_x v_x \quad \textrm{(for the explicit scheme)},\label{eq:mod_exp}\\
    v_t + f v_x &=& \frac{1}{2} h f v_{xx}(1+\lambda) + \frac{1}{2}h \lambda f_x v_x \quad \textrm{(for the implicit scheme)},\label{eq:mod_imp} \textrm{ and} \\
    v_t + f v_x &=& \frac{1}{2} h f v_{xx} \label{eq:mod_semi} \,\quad\quad\quad\quad\quad\quad\quad\quad\quad\textrm{(for the semi-discrete scheme)},
\end{eqnarray}
where $\lambda = \lambda(x,t) = f(x) \frac{k}{h}$ rather than a
constant. The first term on the right-hand side represents the numerical viscosity of each scheme.
For the constant $f$ case, taking $k<\hat{k}$ makes $\lambda<1$,
and decreasing $k$ increases the viscosity in the explicit scheme,
while $(1+\lambda) > (1-\lambda)$ implies higher numerical viscosity in the implicit scheme.
The same argument also holds for a {\em slowly varying} advection speed $f(x)$;
Figure~\ref{fig:imp_exp_win}a shows the errors of both schemes measured in the time slice $T=1.5$ for $f(x)=\frac{1}{2x+1}$.
The errors are plotted for a range of $h$ values, with a CFL prescribed time-step $k=\hat{k}$ used in both schemes,
showing that the implicit scheme has lower accuracy even if we are willing to take as many time-steps as in the explicit case.

However, for {\em strongly varying} advection speeds, this argument is far less convincing.
If $f$ is time-independent
, then \eqref{eq:semi-discr} is a constant coefficient linear system with eigenvalues
$\{ -F_i / h  \, \mid \, i=1,\ldots, M\}$, and $(-\hat{f} / h)$ characterizes its stiffness.
It is natural to expect implicit methods to perform much better on stiff problems, and this turns out to be the case here as well
.  The CFL-prescribed time $\hat{k}$ is still based on $\hat{f}$ even though
the averaged $\bar{f} = \int_0^1 f(x) \, dx$ might be much smaller.  In such scenarios, we might have $\lambda \ll 1$ on most of $\cdomain$, resulting in a relatively large numerical viscosity of explicit scheme on most of $\cdomain.$  Wherever $f \ll \hat{f}$, the additional increase in viscosity due to 
the implicitness of this scheme is relatively modest, and the opportunity to use $k>\hat{k}$ becomes attractive.

\begin{figure}[ht]
    \centering
    \subfloat[]{\includegraphics[width=0.25\columnwidth]{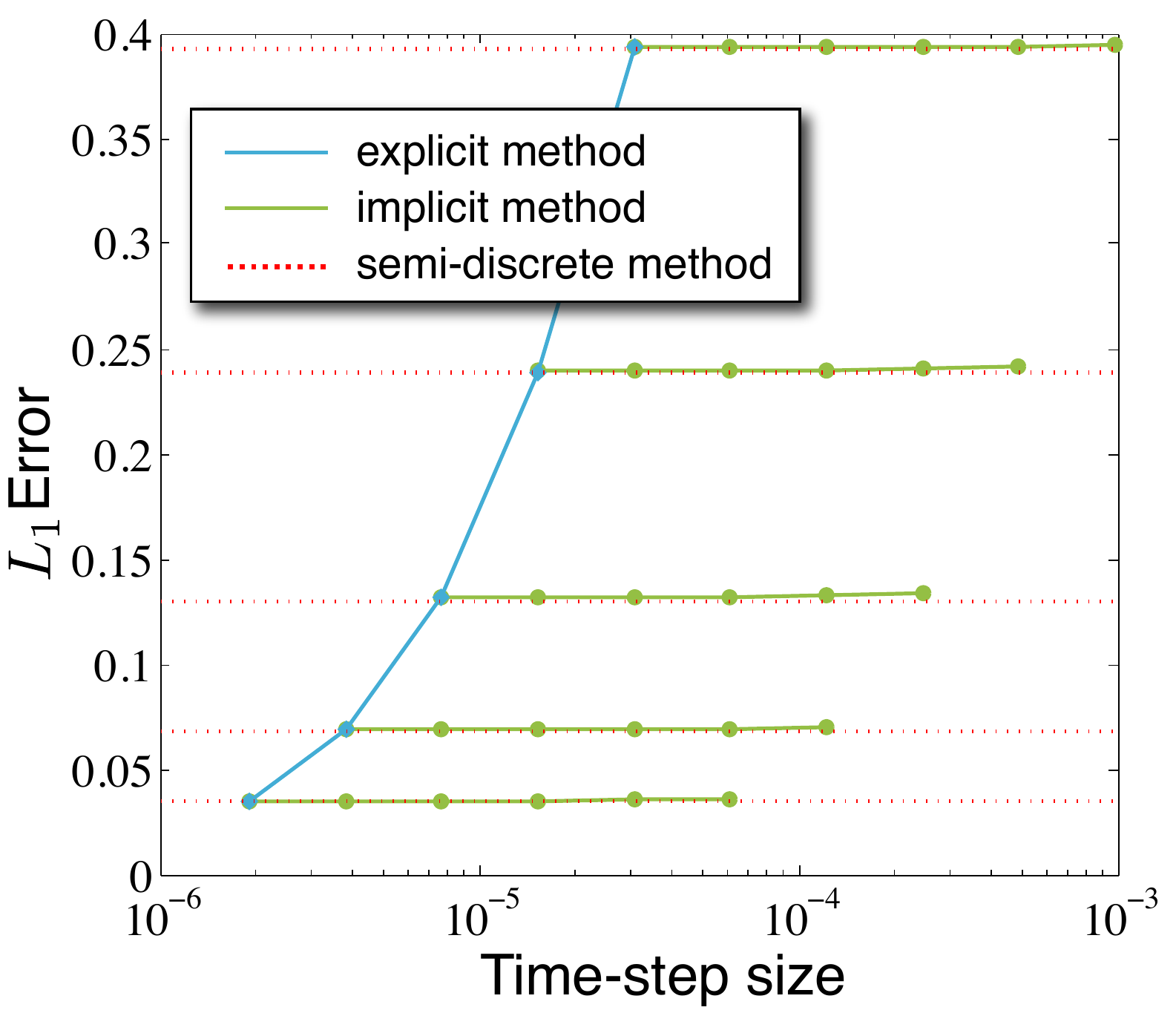}\label{fig:1d_case3}}
    \subfloat[]{\includegraphics[width=0.25\columnwidth]{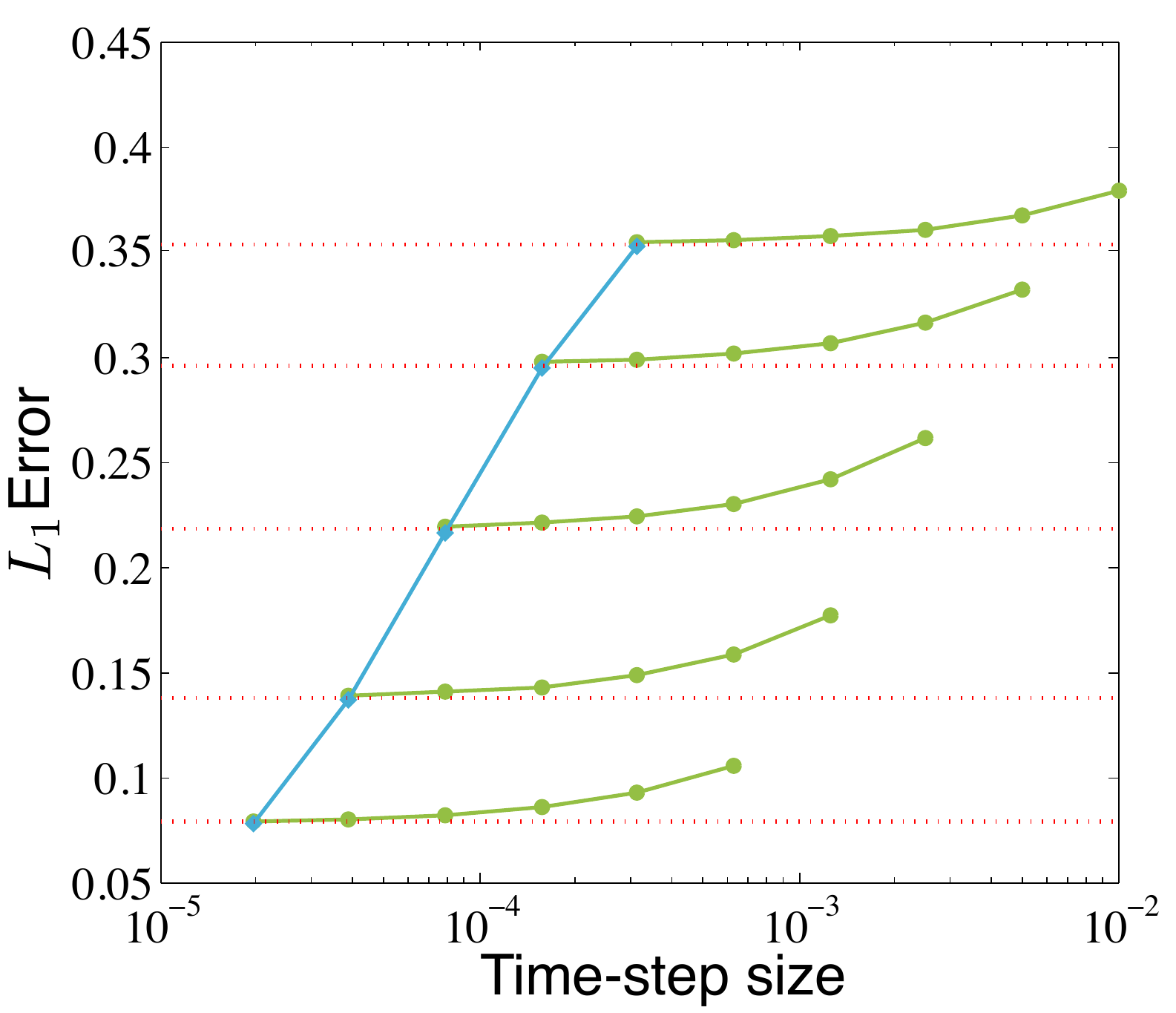}\label{fig:1d_case5}}
    \subfloat[]{\includegraphics[width=0.25\columnwidth]{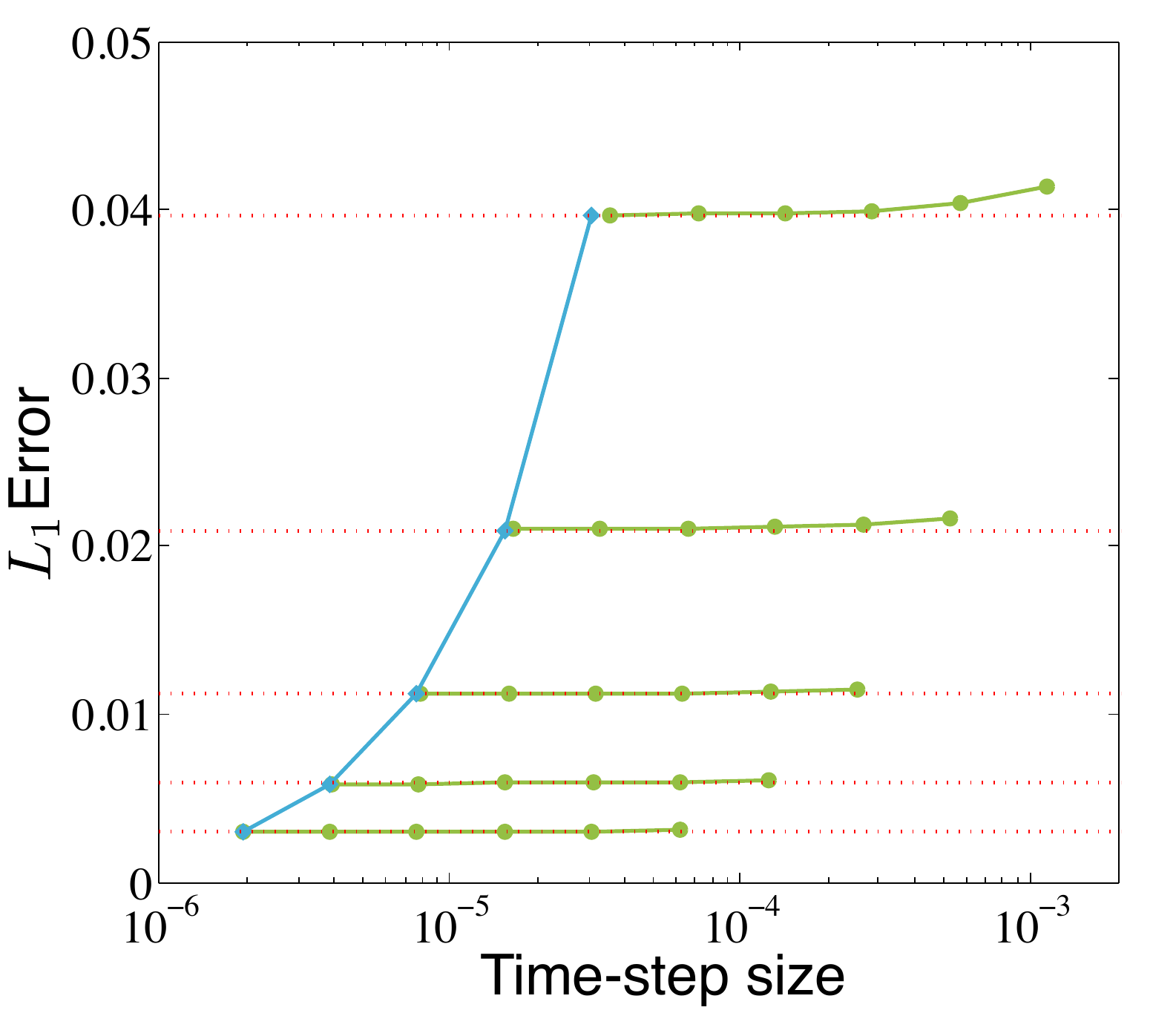}\label{fig:1d_case4-10}}
    \subfloat[]{\includegraphics[width=0.25\columnwidth]{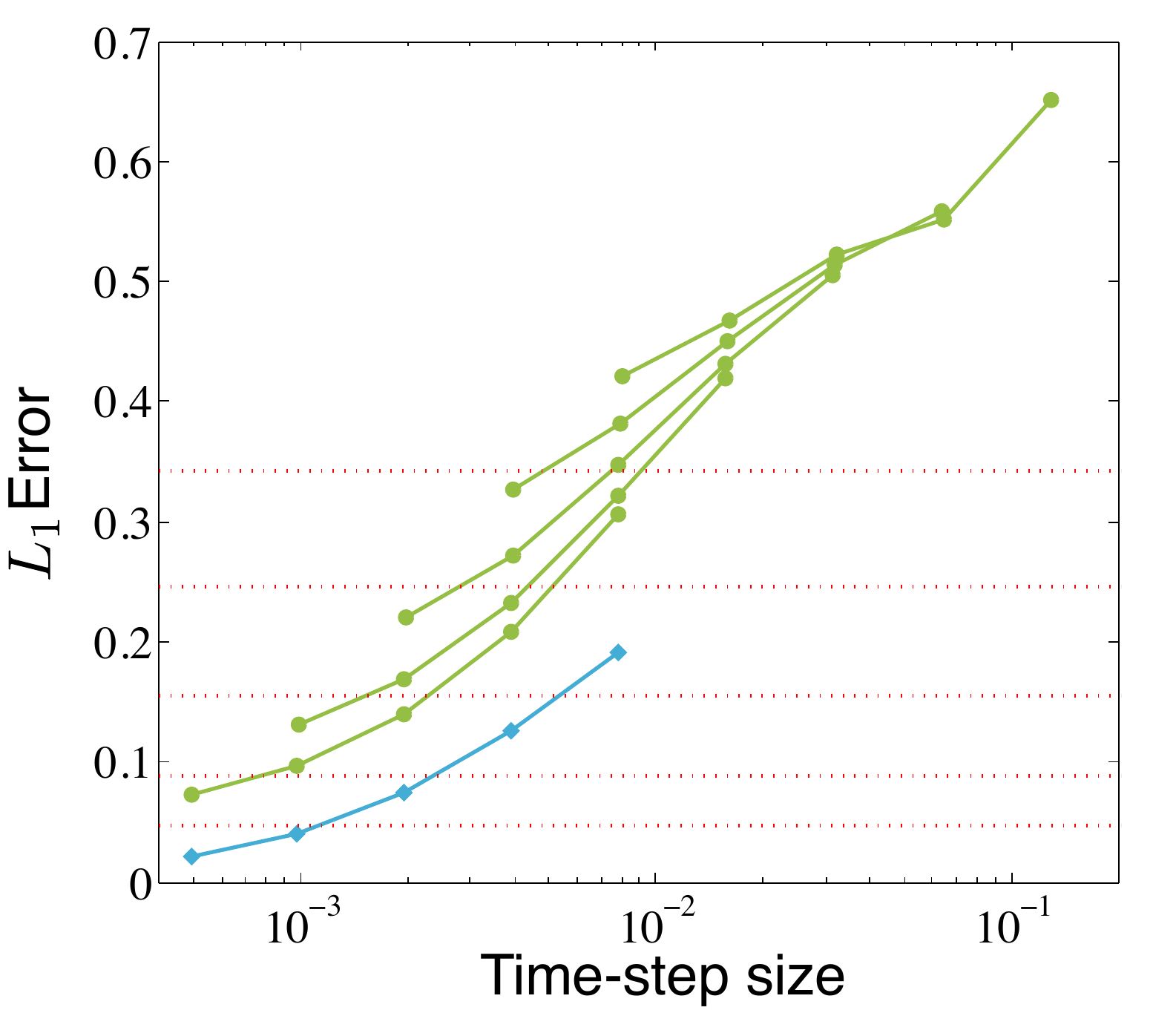}\label{fig:1d_case4-2}}
    \caption{ \textbf{Explicit versus implicit methods:}
    $L_1$ errors at $T=1$
    using different time-step sizes in explicit and implicit schemes.
    The solved equation is $v_t + fv_x = 0$ 
    with the boundary condition
    $v(0,t) = \sin(10\pi t)$. The speed functions: (a) $f=(1+x)^{10}$, (b) $f=100/(50.5+49.5\cos(2\pi x))$, (c) $f=(2-x)^{10}$, and (d)$f=(2-x)^2$.
    Each point on the blue (explicit scheme) curve corresponds to a particular
    $h$ and the largest CFL-prescribed time step $\hat{k}$.  Each green
    (implicit scheme) curve corresponds to a single $h$ but with a number of
    different $k$ values.
    The five green curves from bottom to top correspond to the spatial resolutions of $1024$, $512$, $258$, $128$ and $64$ gridpoints respectively.
    \label{fig:imp_vs_exp}}
\end{figure}

Fig.~\ref{fig:imp_vs_exp} illustrates this observation with several test problems.
All the examples solve the equation $v_t + f v_x = 0$ on $x\in(0,1)$ with the
same boundary condition $v(0,t) = \sin(10\pi t)$ but with different speed functions.
Examples (a-c) use highly stiff speed functions, which strictly
limit the time-step sizes for explicit schemes. Consequently, errors based on spatial
discretization are dominant, and we can take larger time-steps in implicit scheme
without loss of too much accuracy.  In each subfigure a blue curve shows the accuracy of the explicit scheme for various $h$ values and the corresponding CFL-prescribed values of $\hat{k}$.  Each green curve shows the accuracy of the implicit scheme for a fixed $h$, but with a number of different $k\geq \hat{k}$ values.  (Starting from $k=\hat{k}$ and doubling the time-step as we move to the right along each green curve.)  The red dotted lines show the errors of the semi-discrete scheme \eqref{eq:semi-discr} for the same $h$ values.  All errors are measured in $L_1$ norm and in the $T$-time slice only.  (Measuring errors in $L_{\infty}$ norm and/or in all the time slices 
produces very similar results.)

In the first 3 examples (Fig.~\ref{fig:1d_case3}-\ref{fig:1d_case4-10}),
the green curves show that the implicit scheme with time-steps 8$\times$ or even 16$\times$ larger than $\hat{k}$ still yields comparable accuracy.  The ``computational-cost-per-accuracy'' comparison is slightly more subtle since each step to the right along any green curve reduces the cost by a factor of 2, while a step to the right along the blue curve reduces the cost by a factor of 4.  Nevertheless, it is clear that the implicit method is significantly more efficient here.

On the other hand, the last example in Fig.\ref{fig:1d_case4-2} illustrates the performance of both schemes on non-stiff problems. In this case, the implicit scheme exhibits a steep increase in errors as $k$ increases, and the explicit scheme is clearly much better.

\begin{remark}
\label{rem:imp_cost_higher_dim}
The comparison of computational costs based on such plots is also dimension-dependent even if the cost-per-gridpoint-value-update remains the same for explicit and implicit schemes.  In $R^d$ each step to the right along a blue curve would decrease the cost by a factor of $2^{d+1}$, while on green curves it would remain a factor of 2.  This limits the advantages of our proposed techniques in higher dimensions, since the minimum stiffness of the problem needed to justify the use of implicit schemes grows exponentially with $d$.
\end{remark}

\begin{remark}
\label{rem:stiffness_and_more}
Our notion of problem stiffness informally defined in terms of $\hat{f} / \bar{f}$  is a good indicator for the potential usefulness of implicit schemes, but does not fully determine the outcome.  Several additional considerations are enumerated below:
\begin{enumerate}
\item
The location where $\hat{f}$ is attained clearly also matters.  E.g., if it occurs only close to the outflow boundary, it has a smaller impact on the $L_1$ errors.
\item
The accuracy of each method also depends on the rate of change of the boundary condition $\beta(t)$ (and, in a more general case, on $f_t(x,t)$).
In general, if $\beta$ significantly changes on the time scale of $h/\hat{f}$, the stability condition is not really restrictive and the explicit method is likely more efficient.
\item
For a somewhat more general equation,
\begin{equation}
\label{eq:advect_with_source}
v_t + f(x,t) v_x = g(x,t),
\end{equation}
the implicit scheme might produce smaller errors than the explicit even when using the same $k \leq \hat{k}$ and even if the speed function is not really stiff.  See Figure~\ref{fig:imp_exp_win}b and the discussion of semi-Lagrangian methods below.
The numerical methods suitable for the more general stiff semi-linear source term (i.e., the case $g=g(x,t,v)$) are discussed in \cite{LeVequeYee}.
\item
If $\beta$ is constant, both explicit and implicit schemes could be considered as iterative methods to recover
the solution of a stationary problem.  Since the implicit scheme corresponds to Gauss-Seidel iterations, it will 
converge faster.  More generally, it is easy to show that this advantage of implicit schemes also occurs whenever $v(x,t)$ is linear in time, regardless of the amount of stiffness present in the problem.
\end{enumerate}
\end{remark}

\paragraph{Comparison with semi-Lagrangian schemes:}
When interested in large time steps in problems such as \eqref{eq:advect_with_source}, one typical approach is to employ the standard {\em semi-Lagrangian} schemes \cite{FalconeFerretti_book}.  The idea is to follow the characteristic from $(x_i, t_{n+1})$ for time $k$ until reaching some point $(\tilde{x}, t_n)$, where $\tilde{x}$ falls between gridpoints $x_j$ and $x_{j-1}$ and can be recovered as their linear combination  $\tilde{x} = \xi_1 x_j + \xi_2 x_{j-1}$, with $\xi_1, \xi_2 \geq 0$ and $\xi_1 + \xi_2 = 1$.  The first-order semi-Lagrangian scheme is then obtained by approximating $\tilde{x} \approx x_i - k f(x_i, t_{n+1})$ and
\begin{equation}
\label{eq:semiL}
V_i^{n+1} = k g(x_i, t_{n+1}) + \left(\xi_1 V_j^n + \xi_2 V_{j-1}^n \right),
\end{equation}
where the first term in the sum approximates the integral of $g(x,t)$ along the characteristic.
This scheme is also unconditionally stable, but uses an extended stencil (since $i$ and $j$ are generally not the same).
We note that the finite difference schemes considered above can be also interpreted as semi-Lagrangian: the explicit scheme is recovered by simply taking $k \leq \hat{k}$, ensuring that $i = j$;  the implicit scheme corresponds to following the linearized characteristic for a smaller time $\tau_i^{n+1} = kh / [kf(x_i,t_{n+1}) \, + \,h]$ with a subsequent interpolation on the segment between $(x_i,t_n)$ and $(x_{i-1}, t_{n+1})$; see Figure~\ref{fig:imp_exp_win}c.

We note that there are several different sources of errors in semi-Lagrangian techniques:\\
1) due to approximating the integral of $g$;
2) due to approximating the location of $\tilde{x}$; and
3) due to interpolating the value of $V(\tilde{x}, t_n)$.

With $g=0$, only the second and the third of these are present (the scenario illustrated in Figures~\ref{fig:imp_exp_win}a and ~\ref{fig:imp_vs_exp}).
The example in Figure~\ref{fig:imp_exp_win}b was selected to show what happens when the first of the above is the only source of error.
The implicit scheme wins here due to a smaller error in approximating the integral along the characteristic (since $\hat{k} > \tau = \hat{k}/2$).
We have not performed a systematic comparison of accuracy between \eqref{eq:semiL} and \eqref{eq:imp_scheme_1d} for $k > \hat{k}$, but we expect the latter to be similarly more accurate for problems with strongly varying $g(x,t)$.

\begin{figure}[ht]
    \centering
    \subfloat[]{\includegraphics[width=0.4\hsize]{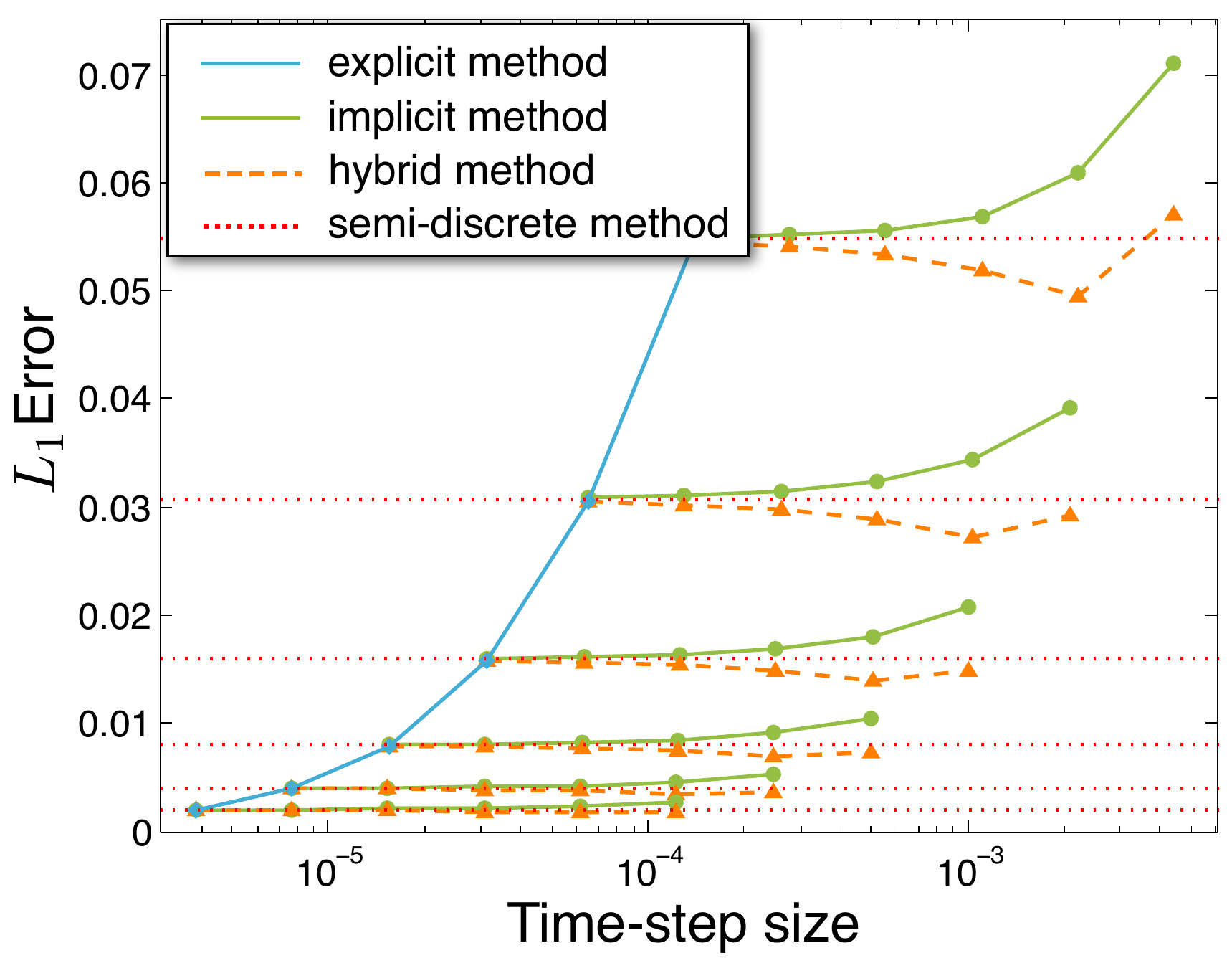}}
    \quad\quad
    \subfloat[]{\includegraphics[width=0.4\hsize]{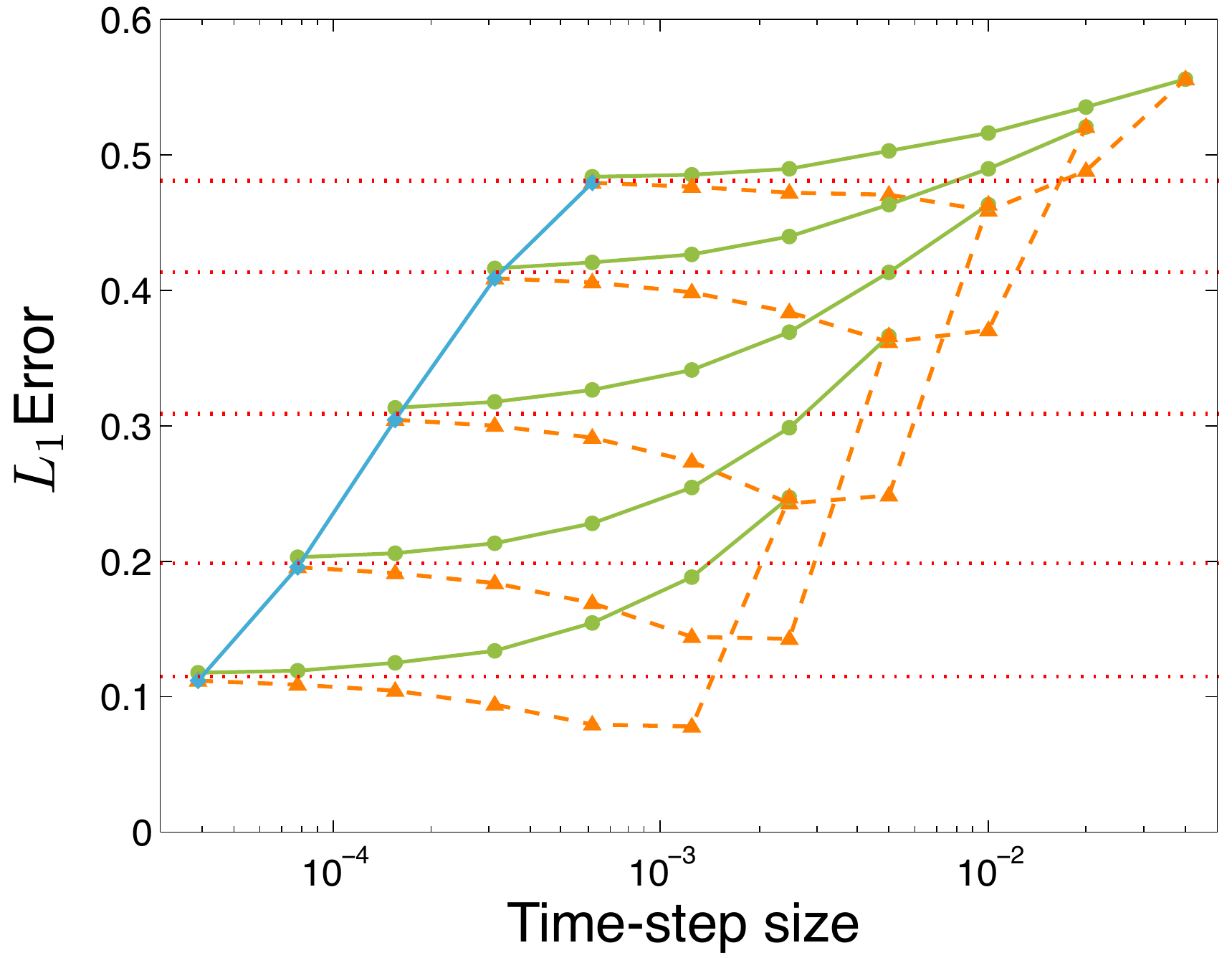}}
    \caption{\textbf{Hybrid method:} $L_1$ error of solving the equation $v_t + f v_x = 0$ with different methods. (a) uses the
    speed function $f=(2-x)^8$ and evaluates the error at $t=0.184$.
    (b) uses the speed function $f=50/(25.5 + 24.5\cos(2\pi x))$ and evaluates the error at $t=0.663$.
    \label{fig:hybrid_vs_imp_vs_exp}}
\end{figure}

\paragraph{Hybrid method:}
We note that it is also possible to combine the explicit and implicit
schemes.  The resulting combination could be best described as {\em ``implicit-on-demand''} or {\em ``opportunistically-explicit''},
but for the sake of brevity we will refer to it as ``hybrid''.
Given a fixed spatial resolution $h$ and a fixed time-step size $k \geq \hat{k}$,
we solve for the values, $V_i^{n+1}$, at the next time-step using two passes.
In the first pass, we compute $V_i^{n+1}$
using explicit scheme~\eqref{eq:exp_scheme_1d} at the nodes where the CFL
stability condition is satisfied (i.e. $\lambda_i^n\le 1$).  In the second
pass, we use implicit scheme~\eqref{eq:imp_scheme_1d} to compute $V_i^{n+1}$ at
the remaining nodes where the CFL stability condition can not be satisfied
(i.e. $\lambda_i^n>1$).
Looking back to the modified equations, this further reduces the part of the
domain on which the numerical viscosity is large.  Let $\domain_e(k)
\subset \cdomain$ be a set on which the CFL would be satisfied with current $h$
and $k$ values.  When $k > \hat{k}$, the hybrid method results in less
numerical viscosity on $\domain_e(k)$ than 
would be produced by the purely
explicit scheme with time-step $\hat{k}$.  This suggests that for many problems
a hybrid method can actually produce {\em smaller} errors than the explicit
scheme even as we lower the cost by using larger time steps.  Figure
\ref{fig:hybrid_vs_imp_vs_exp} illustrates this point for two different speed
functions $f$.  Orange curves, corresponding to the hybrid method start out
with the same error as the explicit scheme when $k=\hat{k}$ and $\domain_e(k) =
\cdomain$.  As $k$ increases, $\domain_e(k)$ shrinks but $\lambda$ becomes
closer to one on that set, thus reducing the errors.  As $k$ continues to grow,
the larger numerical viscosity (of the implicit scheme) on the rest of the
domain becomes dominant; finally, for $k > h / \min_x f(x),$ we have
$\domain_e(k) = \emptyset$ and the hybrid scheme becomes equivalent to the
purely implicit one.

We note that similar hybrid schemes were previously developed in the computational fluid dynamics community for linear advection problems with ``spatially and/or temporally localized stiffness in wave speeds'' \cite{CollinsColellaGlaz} and then later generalized for unstructured meshes in higher dimensions \cite{ORourkeSahota}.  Our version described above uses a different implicit stencil that is more suitable for the purposes of section \ref{s:discretize}.

\vspace*{2mm}
Having shown the advantages of implicit and hybrid methods in the linear context, we now review the basics of optimal control formulation before introducing similar efficient methods for non-linear Hamilton-Jacobi PDEs.

\section{Hamilton-Jacobi PDEs in optimal control.}
\label{s:HJB_intro}
Hamilton-Jacobi equations arise in a variety of applications, including
optimal control, differential games, and modeling of propagating interfaces.
In the context of deterministic optimal control, the Hamiltonian is convex,
and the value function of the process can be recovered by solving
a first-order Hamilton-Jacobi-Bellman PDE \cite{BardiDolcetta}.
Exit-time optimal control problems describe the task of leaving
the domain $\domain \subset \R^d$ as cheaply as possible,
starting from every possible initial configuration of the controlled
system $\x \in \domain$.
With autonomous running cost and dynamics such problems yield
the value function (the minimal cost associated with an optimal control)
$u(\x)$ satisfying a static PDE
\begin{eqnarray}
- \theta u \, + \,
H(\nabla u, \x) = 0, & \text{ on } & \domain \subset \R^d;\\
u = q, & \text{ on } & \boundary,
\label{eq:sHJB}
\end{eqnarray}
where the Hamiltonian $H$ encodes the dependence of $u$
on the running cost and dynamics,
$\theta \geq 0$ is the rate of time-discounting,
while $q$ is
the additional terminal cost charged for crossing $\boundary$.

More specifically, suppose that $\domain$ is open and bounded
and the vehicle's dynamics inside $\cdomain$
is defined by
\begin{eqnarray}
\nonumber
\y^{\prime}(s) &=& \fB(\y(s), \ba(s)),\\
\y (0) &=& \x \in \domain,
\label{eq:auton_dynamics}
\end{eqnarray}
where $\y(s)$ is the system state at the time $s$,
$\ba(s)$ is the currently used control value
(chosen from a compact set $A \subset \R^m$),
$\x$ is the initial system state,
and $\fB: \cdomain \times A \mapsto \R^d$ is the velocity.
The exit time associated with this control is
\begin{equation}
T_{\mathbf{x}, \mathbf{a}} = \min \{ s \in \R_{+,0} | \y(s) \in \boundary \},
\label{eq:exit-time}
\end{equation}
with the convention that $T_{\mathbf{x}, \mathbf{a}} = +\infty$ if $\y(s) \in \domain$ for all $s\geq0$.
The problem description also includes
the running cost $K: \cdomain \times A \mapsto \R$ and
the terminal (non-negative, possibly infinite) cost
$q: \boundary \mapsto \left(\R_{+,0} \cup \{+\infty\} \right)$.
We will assume that\\
$\bullet \,$ $\fB$ and $K$ are Lipschitz-continuous;\\
$\bullet \,$ $0 < K_1 \leq K(\x, \ba) \leq K_2$ for $\forall \x \in \cdomain, \ba \in A$;\\
$\bullet \,$ $q$ is lower semi-continuous and
$\min_{\boundary} q < +\infty$.\\
This allows to specify the total cost of using the control $\ba(\cdot)$ starting from $\x$:
$$
\J(\x, \ba(\cdot)) = \int_0^{T_{\mathbf{x}, \mathbf{a}}}  e^{-\theta s} \, K (\y(s), \ba(s)) \, ds
\, + \, e^{-\theta T_{\mathbf{x}, \mathbf{a}} } \, q(\y(T_{\mathbf{x}, \mathbf{a}})).
$$
The key idea of dynamic programming is to introduce the {\em value function}
\begin{equation}
u(\x) \; = \; \inf_{\ba(\cdot) \in \mathcal{A}} \J(\x, \ba(\cdot)).
\label{eq:value}
\end{equation}
Bellman's optimality principle and Taylor series can be used to
formally derive a static Hamilton-Jacobi-Bellman PDE:
\begin{eqnarray}
\nonumber
- \theta u \, + \,
H(\nabla u, \x) \, = \, - \theta u \, + \,
\min\limits_{\ba \in A}
\left\{
K(\x, \ba) + \nabla u (\x) \cdot \fB ( \x, \ba )
\right\}
\, = \, 0,
&& \text{ for $\x \in \domain$};\\
u(\x) \, = \, q(\x), && \text{ for $\x \in \boundary$}.
\label{eq:HJB_general_static}
\end{eqnarray}
Unfortunately, a smooth solution to Eqn.~\eqref{eq:HJB_general_static}
might not exist even for smooth $\fB$,$K$, $q$, and $\boundary$.
Generally, this equation has infinitely
many weak Lipschitz-continuous solutions, but the unique {\em viscosity solution}
can be defined using additional conditions on smooth test functions \cite{CranLion, CranEvanLion}.
It is a classic result that the viscosity solution of this PDE
coincides with the value function of the above control problem
and the characteristic curves of this PDE coincide with the trajectories of optimal motion;
see \cite{BardiDolcetta} for a detailed discussion and extensive references\footnote{
With $\theta>0$, an optimal control $\ba(\cdot)$ might have $T_{\mathbf{x}, \mathbf{a}}=+\infty$.
Thus, this setting could be more accurately described as ``infinite horizon with time-discounting inside $\domain$
or finite termination upon reaching $\boundary$'', but we stick to the ``exit-time optimal control''  nomenclature
for the sake of brevity.  Our main focus is on the non-time-discounted case (i.e., $\theta = 0$),
where $T_{\mathbf{x}, \mathbf{a}}$ is finite for every optimal control.
}.

For $K \equiv 1$ and $\theta=0$, this control problem amounts to finding the time-optimal trajectories
and $q$ is interpreted as an exit time penalty.
We note that for autonomous dynamics,
the min-time-from-$\x$-to-$\boundary$ problem
is equivalent to the min-time-from-$\boundary$-to-$\x$ problem, with $q$ simply interpreted
as ``entry time penalty'' in the latter case.
However, in the non-autonomous situation where
$\fB = \fB(\y(s), \ba(s), s)$, these problems are quite different.
Both are considered in detail in \cite{VladTimeD}.
While the min-time-from-$\x$-to-$\boundary$ requires a time-dependent
value function and an evolutive PDE, the min-time-from-$\boundary$-to-$\x$
has a stationary value function $u(\x)$ satisfying the static PDE
\begin{eqnarray}
\nonumber
H(\nabla u, \x, u) \, = \, \min\limits_{\ba \in A}
\left\{
1 - \nabla u (\x) \cdot \fB \left( \x, \ba, u(\x) \right)
\right\}
\, = \, 0,
&& \text{ for $\x \in \domain$};\\
u(\x) \, = \, q(\x), && \text{ for $\x \in \boundary$}.
\label{eq:HJB_general_from_boundary}
\end{eqnarray}
In this setting, $s$ is the time since {\em entering} $\cdomain$,
and all optimal trajectories run from $\boundary$ {\em into} the domain,
hence the minus sign in the Hamiltonian.  Along each optimal trajectory,
$s=u(\y(s))$, which explains the third argument of $\fB$ in \eqref{eq:HJB_general_from_boundary}.

In {\em fixed-horizon problems}
the process starts at $\x \in \domain$ at the time $t$ and
stops either at the pre-specified terminal time $T$
or upon reaching the boundary before the time $T$.
Correspondingly, the exit cost $q$ is now charged
on $\boundary \times [0,T] \, \bigcup \, \domain \times \{T\}$.
Assuming no time-discounting (i.e., $\theta=0$) and
defining the value function $v(\x,t)$ similarly to
the exit-time problems, we can again use Bellman's optimality
principle and Taylor series to formally derive the following
evolutive PDE:
\begin{eqnarray}
\label{eq:tHJB}
v_t + H(\nabla v, \x, t) = 0, & \text{ on } & \domain \times [0,T);\\
v = q, & \text{ on } & \boundary \times [0,T] \, \bigcup \, \domain \times \{T\},
\end{eqnarray}
with the Hamiltonian
$$
H(\nabla v, \x, t) \, = \, \min\limits_{\ba \in A}
\left\{
K\left(\x, \ba, t \right) + \nabla v (\x,t) \cdot \fB \left( \x, \ba, t \right)
\right\}.
$$
We refer to \eqref{eq:tHJB} as a {\em terminal-value problem}
since $v(\x,T)$ is known and $v(\x,t)$ needs to be found for earlier
times.  Intuitively, the information is propagating
from the future into the past and one natural approach is to use an explicit
time-marching discretization.  As discussed in section \ref{s:discretize},
this yields a simple computational approach, but results in a restrictive
CFL-stability condition.

We note that for fixed-horizon problems the value function is usually time-dependent
even if $H = H(\nabla v, \x)$ and $q$ is autonomous on $\boundary \times [0,T]$.
The reason is that the remaining time $(T-t)$ might be insufficient to reach
$\boundary$ from $\x$.  One obvious exception is the case $q(\x,T) = + \infty$,
in which $v(\x,t)$ will not depend on time wherever $v$ is finite.  (I.e., all
characteristics with finite values of $v$ reach $\boundary$ before the time $T$.)

\begin{remark}\label{rem:HJB_as_advect}
The above optimal control setting obviously also describes the total cost of ``non-controlled'' processes.
This corresponds to the simplest case where the set $A = \{ \ba \}$ is a singleton,
$\fB(\x,t) = \fB(\x,t, \ba)$ and $K(\x, t) = K(\x, \ba, t)$, resulting in a linear HJB equation.
In 1D and with the running cost $K=0$, this HJB equation \eqref{eq:tHJB} would reduce to a simple non-homogeneous advection equation \eqref{eq:advect} considered in section \ref{s:advection}.
\end{remark}

Several simple types of controlled dynamics are frequently relevant
in many applications.    We will say that the system
has {\em geometric dynamics} if $A = \{ \ba \in \R^d \mid |\ba|=1 \}$ and
$\fB(\x, \ba, t) = f(\x, \ba, t) \ba$, where $\ba$ is the direction of motion,
and $f$ is the speed.
The {\em local controllability} assumption essentially states
that we can move with non-zero speed in any direction;
i.e., $0 < F_1 \leq f(\x, \ba, t) \leq F_2$ for all
$x \in \cdomain, \ba\in A, t \in [0,T].$
For geometric dynamics, ``chattering controls'' can also be used to effectively ``stay in place''.
This ``cost of not moving'' can be defined as
$K(\x,\bm{0},t) = \inf_\mu \left\{ \int_A K(\x,\ba,t) \, d\mu \right\},$
where the infimum is taken over all probability measures $\mu$ on $A$ such that $\int_A f(\x, \ba, t) \ba \, d\mu = \bm{0}$.

If we further assume that the problem is {\em isotropic} (i.e.,
$f = f(\x, t)$ and $K = K(\x, t)$) then the minimization in equation
\eqref{eq:tHJB} can be performed analytically, yielding the
following {\em time-dependent Eikonal equation}:
\begin{equation}
v_t + H (\nabla v, \x, t) = v_t + K(\x,t) -  f(\x, t) | \nabla v | = 0.
\label{eq:tEikonal}
\end{equation}

In autonomous exit-time problems, isotropic cost and dynamics
reduce equation \eqref{eq:HJB_general_static} to a static Eikonal PDE
$K(\x) - f(\x) |\nabla u | = 0$.  It is a simple observation that
this value function can be also interpreted as solving a time-optimal
problem
\begin{equation}
1 - \hat{f}(\x) |\nabla u | = 0,
\label{eq:plain_Eikonal}
\end{equation}
where
$\hat{f}(\x) = f(\x) / K(\x)$.
Finally, for a non-autonomous min-time-from-$\boundary$-to-$\x$
problem \cite{VladTimeD}, the isotropy assumptions reduce equation
\eqref{eq:HJB_general_from_boundary}
to a slightly more general static Eikonal PDE
\begin{equation}
H(\nabla u, \x, u) \, = \, 1 - | \nabla u (\x) | f \left( \x, u(\x) \right)
\, = \, 0.
\label{eq:Eikonal_from_boundary}
\end{equation}


\section{Discretizations and numerical methods.}
\label{s:discretize}

We start by considering two (first order accurate in time)
semi-discretizations of the fixed horizon problem \eqref{eq:tHJB}
with geometric dynamics
(i.e., $A= S_1$ and $\fB \left( \x, \ba, t \right) = f ( \x, \ba, t) \ba$).
We will assume that the time step is $k=T/N$ and a typical time-slice is
$t_n = nk \leq t_N = Nk = T.$  We will consider a sequence of functions $v^n:\cdomain \mapsto R$
approximating $v(\x, t_n)$.  For notational simplicity we will use the superscripts to specify the
time slice in the cost and dynamics (e.g., $K^n(\x, \ba) = K(\x, \ba, t_n)$).
Since \eqref{eq:tHJB} is a terminal time problem (i.e., $v^N(\x) = q^N(\x)$ is specified on $\cdomain$),
it is logical to consider the problem of finding $v^n(\x)$ when $v^{n+1}(\x)$ is already known.
Both time-discretizations are obtained by replacing
$v_t(\x,t_n)$ with the divided difference $(v^{n+1}(\x)-v^n(\x))/k$.
The explicit formula results from evaluating the Hamiltonian using $v^{n+1}(\x)$:
\begin{equation}
\label{eq:exp_semi}
v^n \; = \; v^{n+1} \, + \, k H \left( \nabla v^{n+1}, \x, t_{n} \right)
.
\end{equation}
Similarly, the implicit formula results from evaluating the Hamiltonian using $v^{n}(\x)$:
\begin{equation}
\label{eq:imp_semi}
v^n  = v^{n+1} \, + \, k H \left( \nabla v^{n}, \x, t_{n} \right)
\; = \; v^{n+1} \, + \, k
\min\limits_{|\ba|=1}
\left\{
K^n\left(\x, \ba \right) + \left( \nabla v^{n} \cdot \ba \right) f^n \left( \x, \ba \right)
\right\}.
\end{equation}
The key idea of our approach is to reinterpret \eqref{eq:imp_semi} as a boundary value problem
for $v^n(\x)$ that corresponds to an auxiliary exit-time optimal control problem on $\domain$.
In particular, \eqref{eq:imp_semi} can be rewritten as \eqref{eq:HJB_general_static} if we define
$$
\fB ( \x, \ba ) = f ( \x, \ba ) \ba, \qquad
\theta = \frac{1}{k}, \qquad
K \left(\x, \ba \right) = K^n\left(\x, \ba \right) + \frac{v^{n+1}(\x)}{k}, \qquad
\text{and } q(\x) = q^n(\x).
$$

Another ``stationary reinterpretation'' is even more convenient for the case when the running cost $K^n$ is isotropic.
If the optimal strategy at $(\x, t_n)$  is to stay in place,
\eqref{eq:imp_semi} is equivalent to $v^n(\x) =  v^{n+1}(\x) \, + \, k K^n\left(\x \right)$.
But on the rest of $\domain \times \{t_n\}$, we have $v^{n+1}  - v^n + k K^n > 0$,
and $v^n$ can be reinterpreted as a a value function
for an auxiliary ``min-time-from-$\boundary$-to-$\x$'' optimal control problem.
Indeed,
\begin{eqnarray*}
\frac{v^n(\x) - v^{n+1}(\x) - k K^n(\x)}{k} &=&
\min\limits_{|\ba|=1}
\left\{
\left( \nabla v^{n}(\x) \cdot \ba \right) f^n \left( \x, \ba \right)
\right\}
\qquad \Longrightarrow \\
0 &=& \min\limits_{|\ba|=1} \left\{
1 + \left( \nabla v^{n}(\x) \cdot \ba \right)
\frac{k f^n \left( \x, \ba \right)}{ v^{n+1}(\x)  - v^n(\x) + k K^n(\x)}
\right\}.
\end{eqnarray*}
The latter is clearly equivalent to \eqref{eq:HJB_general_from_boundary} if we define
the modified velocity
$$
\fB \left( \x, \ba, s \right) = \tilde{f} ( \x, -\ba, s) \ba,
\qquad 
\qquad
\qquad
\tilde{f} ( \x, \ba, s) \, = \,
\frac{k f^n \left( \x, \ba \right)}
{v^{n+1}(\x) \, - \, s \, + \, k K^n(\x)},
$$
where $s$ is the time since entering $\cdomain$.  In the next subsection this idea is used to build a grid discretization for the case of isotropic dynamics.

For any grid discretization of the explicit equation \eqref{eq:exp_semi},
if the local stencil is used to approximate $\nabla v^{n+1}$, this will yield CFL-type stability
conditions, restricting the time step and thus increasing the computational cost.
In contrast, any consistent grid discretization of the implicit equation \eqref{eq:imp_semi}
will be unconditionally stable.
In addition, it is also possible to employ a {\em hybrid} approach by using  \eqref{eq:exp_semi}
on some part of the domain $\domain_e \subset \domain$ (wherever the CFL condition is satisfied),
then taking the result to specify additional boundary conditions on $\boundary_e$,
and using \eqref{eq:imp_semi}
to define $v^n$ on $\domain \backslash \domain_e$.
The resulting method (fully discretized on  a Cartesian grid $X$) is summarized in
Algorithm \ref{alg:imp_and_hybrid}.
The opportunity to use larger time steps is the key advantage of implicit and hybrid methods,
but the overall efficiency obviously hinges on our ability to quickly solve
the boundary value problem (either \eqref{eq:HJB_general_static} or \eqref{eq:HJB_general_from_boundary}).
As we explain in subsection \ref{ss:labelset}, fast non-iterative methods for the latter problem are available,
but in the general anisotropic case they rely on extending the discretization stencil,
which makes it harder to compare the performance/efficiency with the time-explicit approach.
To simplify this comparison, we will focus on isotropic problems, for which non-iterative methods are applicable
even when \eqref{eq:imp_semi} is discretized on a local stencil.

\subsection{Eulerian discretizations for the isotropic case}

We will consider a uniform Cartesian grid $X$ superimposed on
the domain $\domain \subset R^d$.
For the sake of notational simplicity, we will assume that $d=2$ though
generalizations to nonuniform and higher dimensional grids is straightforward.
If $h$ is the gridsize,
a typical gridpoint will have coordinates $\x_{ij} = (x_i, y_j) = (ih, jh)$,
and the total number of gridpoints in $X$ is $M = O(h^{-2})$.
We will also assume that the domain boundary is conveniently discretized by the grid
(e.g., all examples in section \ref{s:experiments} are considered on a
grid-aligned rectangular domain $\domain$).

We will use a grid function $U_{ij} = U(\x_{ij})$ to approximate
the viscosity solution $u(\x_{ij})$ of the static PDE \eqref{eq:Eikonal_from_boundary}.
Employing the usual one-sided
approximations of partial derivatives
$$
u_x(x_i, y_j) \approx D^{\pm x}_{ij} U =
\frac{ U_{i \pm 1, j} - U_{i,j} }{ \pm h}; \qquad
u_y(x_i, y_j) \approx D^{\pm y}_{ij} U =
\frac{ U_{i, j \pm 1} - U_{i,j} }{ \pm h},
$$
the standard upwind discretization of \eqref{eq:Eikonal_from_boundary}
can be written as follows \cite{SethBook2}:
\begin{equation}
f\left( \x_{ij}, U_{ij} \right) 
\left[
\begin{array}{l}
\left(\max \left( D^{-x}_{ij}U, \, -D^{+x}_{ij}U, \, 0 \right)\right)^2 \\
\, + \, \left(\max \left( D^{-y}_{ij}U, \, -D^{+y}_{ij}U, \, 0 \right)\right)^2
\end{array}
\right]^{1/2}
\; = \; 1.
\label{eq:upwind_discr_Eikonal_from_boundary}
\end{equation}
This discretization is consistent and monotone,
which can be used to prove the convergence of $U$ to $u$;
see \cite{BarlesSouganidis}.
If $f$ depends only on its first argument,
\eqref{eq:upwind_discr_Eikonal_from_boundary}
reduces to
a (quadrant-by-quadrant) quadratic equation for $U_{ij}$,
which can be efficiently solved by a Fast Marching Method in $O(M \log M)$ operations;
see \cite{SethFastMarcLeveSet, SethBook2, SethSIAM}.
A modified version of this method is described in section \ref{ss:algorithms}.

For time-dependent problems, we
will use a grid function
$V^n_{ij} = V(\x_{ij}, t_n)$ to approximate
the viscosity solution $v(\x_{ij}, t_n)$ of PDE \eqref{eq:tEikonal},
where $t_n = nk \leq t_N = Nk = T.$
Since this is a terminal value problem,
and the information propagates backward in time, we will
need to compute all $V^n_{ij}$ values with all $V^{n+1}_{ij}$ values already known.
We will approximate the time derivative with the usual first-order divided difference
$v_t(x_i, y_j, t_n) \approx (V^{n+1}_{ij} - V^n_{ij}) / k.$
The fully discrete version of the explicit scheme \eqref{eq:exp_semi} can be then written as
\begin{equation}
\frac{V_{ij}^{n+1} - V_{ij}^{n}}{k} \, + \, K(\x_{ij}, t_n)
\, - \, f \left( \x_{ij}, t_n \right)
\left[
\begin{array}{l}
\left(\max \left( D^{-x}_{ij}V^{n+1}, \, -D^{+x}_{ij}V^{n+1}, \, 0 \right)\right)^2\\
\, + \, \left(\max \left( D^{-y}_{ij}V^{n+1}, \, -D^{+y}_{ij}V^{n+1}, \, 0 \right)\right)^2
\end{array}
\right]^{1/2}
\; = \; 0.
\label{eq:explicit_discr}
\end{equation}
This is a linear equation for $V_{ij}^n$, and the explicit causality of the method
results in the computational cost of $O(M)$ per time slice.
The convergence to viscosity solution is again demonstrated
by an argument in \cite{BarlesSouganidis}.
However, this method is only conditionally stable, since the monotonicity of
\eqref{eq:explicit_discr} can be guaranteed only if $k F_2 \leq h / \sqrt{2}$.
Thus, to compute the solution for $t \in [0, T]$, the total cost is $O(\frac{M}{h} \sqrt{2} F_2)$
even though the actual speed $f$ might be much smaller than $F_2$ on most of
$\domain$.

The reduction of this computational cost is the primary
motivation for our proposed approach.
The fully discrete version of the implicit scheme \eqref{eq:imp_semi} can be then written as
\begin{equation}
\frac{V_{ij}^{n+1} - V_{ij}^{n}}{k} \, + \, K(\x_{ij}, t_n)
\, - \, f \left( \x_{ij}, t_n \right) 
\left[
\begin{array}{l}
\left(\max \left( D^{-x}_{ij}V^{n}, \, -D^{+x}_{ij}V^{n}, \, 0 \right)\right)^2\\
 \, + \, \left(\max \left( D^{-y}_{ij}V^{n}, \, -D^{+y}_{ij}V^{n}, \, 0 \right)\right)^2
\end{array}
\right]^{1/2}
\; = \; 0.
\label{eq:implicit_discr}
\end{equation}
The resulting system of equations is unconditionally stable -- the scheme is monotone
for an arbitrary time-step $k$; see \cite{Souganidis85}.  As before, this equation has to
be solved for each $\x_{ij}$.  But in contrast to \eqref{eq:explicit_discr},
the equation for $V_{ij}^{n}$ is non-linear and the system is coupled.
Our key idea is to treat each time-slice of this system as
a stationary boundary value problem, which is then efficiently solved by
the modified Fast Marching Method; see Algorithm 2.
Note that,
if $\max \left( D^{-x}_{ij}V^{n}, \, -D^{+x}_{ij}V^{n} \right) < 0$
and $\max \left( D^{-y}_{ij}V^{n}, \, -D^{+y}_{ij}V^{n} \right) < 0$, then
\eqref{eq:implicit_discr} reduces to
\begin{equation}
V_{ij}^{n} \; = \; V_{ij}^{n+1} \, + \, k K(\x_{ij}, t_n).
\label{eq:implicit_stay_in_place}
\end{equation}
This corresponds to situations were the optimal control is to ``stay in place'',
making it easy to compute $V_{ij}^{n}$ without knowing the adjacent grid values in the time slice $t_n$.
In all other cases, we can re-write \eqref{eq:implicit_discr}
as follows:
\begin{equation}
1\, - \, \frac{k f \left( \x_{ij}, t_n \right)}
{(V_{ij}^{n+1} - V_{ij}^{n}) \, + \, k K(\x_{ij}, t_n)}
\,
\left[
\begin{array}{l}
\left(\max \left( D^{-x}_{ij}V^{n}, \, -D^{+x}_{ij}V^{n}, \, 0 \right)\right)^2\\
 \, + \, \left(\max \left( D^{-y}_{ij}V^{n}, \, -D^{+y}_{ij}V^{n}, \, 0 \right)\right)^2
\end{array}
\right]^{1/2}
\; = \; 0.
\label{eq:implicit_discr_static}
\end{equation}
We can then consider $V^n$ as the solution of the system
\eqref{eq:upwind_discr_Eikonal_from_boundary}, with the new speed function
\begin{equation}
\label{eq:modified_speed}
\tilde{f} (\x_{ij}, U_{ij}) \; = \; \frac{k f \left( \x_{ij}, t_n \right)}
{(V_{ij}^{n+1} - U_{ij}) \, + \, k K(\x_{ij}, t_n)}.
\end{equation}
This results in a computational cost of $O(M \log M)$ per time-slice,
with the size of time-step dictated by the accuracy considerations alone.
In section \ref{s:experiments} we show that
a relatively large $k$ already results in sufficient precision
for many strongly inhomogeneous problems.  This often makes the method
advantageous despite its higher computational cost per time-slice.





\subsection{Single-pass methods for static HJB equations}
\label{ss:labelset}
Static HJB equations arise in a wide range of applications.  As a result, the efficient numerical methods for them have been an active research area for the last 15-20 years.
The efficiency here is defined in terms of the total number of floating point operations needed to obtain a numerical solution on a fixed grid.  This is somewhat different from the traditional focus of numerical analysis on the rate of convergence of numerical solutions to the viscosity solutions under the grid refinement.
The related literature is rather broad; here we provide only a brief description of the main approaches, with the context of the previous section in mind.  More comprehensive overviews and discussions of the many connections to efficient algorithms on graphs can be found in \cite{ChacVlad,VladMSSP}.

The computational challenge in solving static PDEs stems from the fact that the systems of discretized equations are typically nonlinear and coupled.
Let $N_{ij} = N(\x_{ij})$ be the set of gridpoints adjacent to $\x_{ij}$ and
$NU_{ij} = \{ U(\x) \mid  \x \in N_{ij} \}$ be the set of adjacent gridpoint values.
If $NU_{ij}$ values were already known, the discretized equation (e.g., the equation \eqref{eq:upwind_discr_Eikonal_from_boundary} for the Eikonal equation)
could be used to solve for the value of $U_{ij}$.
Assuming that this (implicitly defined) solution is unique, we will denote it as $U_{ij} = \mathcal{G}_{ij} (NU_{ij})$ and will
say that the vector of grid values $U$ is a fixed point of the operator $\mathcal{G}: \R^M \to R^M$.
Since the $NU_{ij}$ values are a priori unknown, one natural approach is to proceed iteratively:
start with a suitable guess $U^0 \in \R^M$ of all gridpoint values and define $U^{r+1} = \mathcal{G}(U^r).$
Assuming that $\mathcal{G}$ has a unique fixed point, we will say that this process converges in $\rho$ iterations if $U = U^{\rho} = \mathcal{G}(U^{\rho})$
regardless of $U^0$.  The total computational cost on this fixed grid is then $O(\rho M)$.
For many discretizations of static PDEs, if the process is performed in an infinite-precision arithmetic, $\rho$ could be in fact infinite even if these iterations converge.
On a realistic finite-precision computer, this typically translates to a finite $\rho(M)$, which grows with $M$.
We will say that a numerical algorithm is ``single-pass'' (or ``non-iterative'') if there exists some a priori known upper bound on $\rho$ independent of $M$ and of the machine precision.

It is fairly straightforward to show the uniqueness of $\mathcal{G}$'s fixed point and the convergence of the iterative process for a variety of Eulerian and semi-Lagrangian discretizations of \eqref{eq:HJB_general_from_boundary}; e.g., see \cite{Souganidis85, Falcone_InfHor, FalconeFerretti}.  For the first-order upwind discretization of the Eikonal equation \eqref{eq:plain_Eikonal}, it is easy to show that the convergence will be attained after at most $\rho = M$ iterations.  The proof relies on the {\em causal property} of the Eikonal also inherited by this particular discretization: $U_{ij}$ depends only on a subset of {\em smaller} values in $NU_{ij}$; as a result, at least one new gridpoint receives the correct/final value after each iteration.  Gauss-Seidel iterations can be used to significantly decrease the number of iterations-to-convergence, but $\rho$ will then be strongly dependent on the ordering imposed on the gridpoints.  One popular (``Fast Sweeping'') approach
is to conduct the Gauss-Seidel iterations, alternating through the list of geometric gridpoint orderings \cite{Danielsson, BoueDupuis, Zhao, TsaiChengOsherZhao}.  This idea works particularly well when the characteristics rarely change their direction (e.g., when all of them are straight lines) and the computational domain geometry is relatively simple.  For the first order upwind discretization of the Eikonal PDE, the number of needed iterations becomes largely independent from the gridsize $h$ and the overall cost is asymptotically $O(\nu M)$, where $\nu$ is the upper bound on the number quadrant-to-quadrant switches of a characteristic.  Unfortunately, $\nu$ is typically unknown a priori and strongly depends on the grid orientation.  Moreover, for the general anisotropic HJB, where $\rho$ is typically infinite, theoretical bounds on this ``sweeping'' computational cost remain a challenge (and thus such algorithms are not necessarily ``single-pass'' according to the above definition). Still, the experimental evidence shows that Fast Sweeping results in substantial computational savings for a variety of discretizations \cite{Li_1, Li_2} and even for non-convex Hamiltonians \cite{KaoOsherQian}.

An alternative approach is to exploit the causal properties of the discretization to effectively decouple the system and solve equations one at a time.
The classical Dijkstra's algorithm \cite{Diks} uses this approach to solve a related discrete problem of finding the shortest path on a graph with non-negative transition costs.  The key idea is to subdivide the graph nodes into 3 classes: $Accepted$ (for which the exact value is already known), $Considered$ (for which the current/tentative values are available based on their $Accepted$ neighbors only), and $Far$ nodes, for which no tentative value can be reliably computed at that stage of the algorithm.  The causal nature of the problem guarantees that the smallest of the tentative $Considered$ values is actually correct.  Therefore, a typical stage of the algorithm consists of making the corresponding node $Accepted$ and updating the nodes of its not-yet-$Accepted$ neighbors.  For graphs with small/bounded node-connectivity, the overall cost of this algorithm is $O(M \log M)$, where the $\log M$ term results from the need to maintain the min-heap list of $Considered$ nodes. Two Dikstra-like methods were introduced for first-order upwind discretizations of the Eikonal equation on a grid by Tsitsiklis \cite{Tsitsiklis_conference, Tsitsiklis} and Sethian \cite{SethFastMarcLeveSet, SethSIAM}; see \cite{SethVlad3} for a detailed discussion of similarities and differences in these two approaches. Sethian's Fast Marching Method was later extended to higher order accurate schemes \cite{SethSIAM}, to triangulated meshes in $\R^d$ and on manifolds \cite{KimmSethTria, SethVlad1}, and to quasi-variational inequalities \cite{SethVladHybrid}.  First-order upwind discretizations of the general HJB equations are generally not causal; Ordered Upwind Methods \cite{SethVlad2, SethVlad3, AltonMitchell2} dynamically extend the stencil just enough to ensure the causality, enabling space-marching solvers for these more general equations.  In \cite{VladTimeD} Ordered Upwind Methods were extended to  $u$-dependent Hamiltonians \eqref{eq:HJB_general_from_boundary}; applications of Fast Marching to similarly extended isotropic case \eqref{eq:Eikonal_from_boundary} were previously illustrated in \cite{SethVladHybrid}.  Two recent related methods \cite{AltonMitchell2, Mirebeau3} perform the upwind stencil extension as a separate pre-processing step, with the goal of improving both the accuracy and efficiency.

\subsection{Implicit and hybrid numerical methods for time-dependent HJB}
\label{ss:algorithms}
Algorithm \ref{alg:imp_and_hybrid} summarizes our implicit and hybrid solvers for equation \eqref{eq:tHJB}.
Within the main loop of the algorithm, we need to solve an auxiliary static problem \eqref{eq:HJB_general_from_boundary} in each time slice.
For the isotropic case, our implementation accomplishes the latter by a modified version of Fast Marching Method, described in Algorithm \ref{alg:FMM}.
We emphasize that, at least in principle, this could be done by {\em any} of the `` fast'' methods discussed above.
The resulting efficiency would yet again depend on the advantages or disadvantages of a particular solver for this type of static problems. The issue of performance comparison of Marching and Sweeping remains contentious, with each approach advantageous on a separate subclass of examples; see \cite{GremaudKuster, HysingTurek} for experimental comparison and \cite{ChacVlad} for the new (two-scale) methods combining the advantages of both approaches.

\begin{algorithm}[h]
\caption{Implicit (and Hybrid) Methods for time-dependent HJB equations.}
\label{alg:imp_and_hybrid}
\algsetup{indent=2em}
\begin{algorithmic}[1]
\STATE INITIALIZATION:
\STATE $V^N_{ij} := q(\x_{ij}, T)$ for all $(i,j)$
\STATE $n := N$
\STATE
\STATE MAIN LOOP:
\WHILE{$n > 0$}
    \STATE $n := n-1$
    \STATE Set up the non-updated set $Q = \{ \x_{ij} \mid \x_{ij} \in \boundary \}$
           with $V_{ij}^n := q(\x_{ij}, t_n)$.
    \STATE \hspace{\algorithmicindent}
    \vspace*{1mm}
    \hrule
    \vspace*{1mm}
    \STATE
    \STATE HYBRID VERSION ONLY:
	\FOR{each $\x_{ij} \not \in Q$ such that $\quad \max\limits_{\ba \in A} f(\x_{ij}, \ba, t_n) \; \leq \; \frac{h}{k\sqrt{2}} \;$}
        \STATE Compute $V_{ij}^n$ by an explicit update.  {\em (For the Eikonal case: formula \eqref{eq:explicit_discr}.)}
        \STATE Add $\x_{ij}$ to $Q$.
    \ENDFOR
    \STATE \hspace{\algorithmicindent}
    \vspace*{1mm}
    \hrule
    \vspace*{1mm}
    \STATE
    \STATE Solve an auxiliary static problem \eqref{eq:HJB_general_from_boundary} by a ``single pass'' method
    \STATE to recover $V_{ij}^n$ on $X \backslash Q$ treating $Q$ as a computational boundary.
    \STATE {\em (For the Eikonal case: use Fast Marching to solve \eqref{eq:upwind_discr_Eikonal_from_boundary} on $X \backslash Q$}
    \STATE {\em with the speed defined by \eqref{eq:modified_speed}; see Algorithm \ref{alg:FMM}.)}
\ENDWHILE
\end{algorithmic}
\end{algorithm}


\begin{algorithm}[hhhh]
\caption{Fast Marching Method pseudocode for a time slice $t_n$.}
\label{alg:FMM}
\algsetup{indent=2em}
\begin{algorithmic}[1]
\STATE Initialization:
\FOR{each gridpoint $\x_{ij} \in X$}
	\IF{$\x_{ij} \in Q$}
        \STATE Mark $\x_{ij}$ as $Considered$ and add it to the Considered List $L$.
        \STATE (Note that its $V_{ij}^n$ has been already set in Algorithm \ref{alg:imp_and_hybrid}.)
	\ELSE
		\STATE Mark $\x_{ij}$ as $Far$ and set $V_{ij}^n := V_{ij}^{n+1} \, + \, k K(\x_{ij}, t_n)$.
	\ENDIF
\ENDFOR
\STATE End Initialization
\STATE
\WHILE{$L$ is nonempty}
	\STATE Remove the point $\xbar$ with the smallest value from $L$.
    \STATE Mark $\xbar$ as $Accepted$.
	\FOR{ each not-yet-$Accepted$ $\x_{ij} \in N(\xbar) \backslash Q$}
    \IF{$V^n(\xbar) < V_{ij}^n$}
       	\STATE Compute a temporary value $\widetilde{V}_{ij}^n$ using the upwinding discretization \eqref{eq:implicit_discr_static}. \label{algFMM:line:solve}
    \ENDIF
	\IF{$\widetilde{V}_{ij}^n < V_{ij}^n$}
		\STATE $V_{ij}^n := \widetilde{V}_{ij}^n$
	\ENDIF
	\IF{$\x_{ij}$ is $Far$}
		\STATE Mark $\x_{ij}$ as $Considered$ and add it to $L$.
	\ENDIF
	\ENDFOR
\ENDWHILE
\end{algorithmic}
\end{algorithm}

\pagebreak
Line \ref{algFMM:line:solve} of Algorithm \ref{alg:FMM} requires solving the discretized equation \eqref{eq:implicit_discr_static} at a single gridpoint.
It is easy to show that this can be accomplished in a quadrant-by-quadrant fashion and only using $\x_{i,j}$'s $Accepted$ neighbors \cite{SethSIAM}.
Moreover, only quadrants adjacent to the newly accepted gridpoint $\xbar$ are relevant.
E.g., assuming that $\xbar = \x_{i-1,j}$ and using the notation $W_1 = V^n(\xbar)$, $W_0 = V_{ij}^{n+1}$ we consider 4 variants:\\
$\bullet \, $ if both $\x_{i,j+1}$ and $\x_{i,j-1}$ are not-yet-$Accepted$\\
(corresponding to the case $D^{-x} \geq \max(-D^{+x},0)$ and $0 \geq \max(D^{-y},-D^{+y})$
in equation \eqref{eq:implicit_discr_static}),\\
set $\widetilde{V}_{ij}^n$
to be the smallest solution $V$ of the quadratic equation
\begin{equation}
\label{eq:onesided}
\left[ \frac{V-W_1}{h} \right]^2
\; = \;
\left[ \frac{(W_0-V) \, + \, k K(\x_{ij}, t_n)}{k f \left( \x_{ij}, t_n \right)} \right]^2
\end{equation}
satisfying $V \geq W_1.$
If no such solution exists, set $\widetilde{V}_{ij}^n = +\infty$, ensuring that the update will be rejected.
\\ \\
$\bullet \, $ if $\x_{i,j-1}$ is $Accepted$, but $\x_{i,j+1}$ is not, set $W_2 = V_{i,j-1}^n$;\\
$\bullet \, $ if $\x_{i,j+1}$ is $Accepted$, but $\x_{i,j-1}$ is not, set $W_2 = V_{i,j+1}^n$;\\
$\bullet \, $ if both $\x_{i,j+1}$ and $\x_{i,j-1}$ are $Accepted$, set $W_2 = \min \left(V_{i,j-1}^n, \, V_{i,j+1}^n \right)$;\\
Define $\widetilde{V}_{ij}^n$ to be the smallest solution  $V$ of
\begin{equation}
\label{eq:twosided}
\left[ \frac{V-W_1}{h} \right]^2 \, + \, \left[ \frac{V-W_2}{h} \right]^2
\; = \;
\left[ \frac{(W_0-V) \, + \, k K(\x_{ij}, t_n)}{k f \left( \x_{ij}, t_n \right)} \right]^2
\end{equation}
satisfying $V \geq \max(W_1,W_2).$  If no such solution exists, default to formula \eqref{eq:onesided}.

\newcommand{\xx}{\bm{x}}
\newcommand{\HJB}{Hamilton-Jacobi-Bellman}

\section{Numerical Experiments} \label{s:experiments}
In this section, we conduct a number of numerical experiments and compare
the accuracy and performance of the explicit, implicit and hybrid methods.
We consider four time-dependent HJB PDE examples
on a 2D spatial domain $\Omega = [0,1]\times[0,1]$. The numerical solver
is advanced backwards from $t=T>0$ till $t=0$. The time-step size of explicit method
is always chosen based on the CFL condition, while implicit and hybrid methods are tested with
a range of larger time-step sizes.
For each example, we report the numerical error in the $t=0$ time-slice.
The computational time is reported by averaging over $8$ runs for each
algorithm and each parameter configuration. All the timings were measured on a
dual Intel Xeon X5570 (2.93 GHz) processor machine with $48$G memory. The code
\footnote{ The code can be downloaded from \url{http://www.cs.columbia.edu/~cxz/TimeDepHJB/}.}
was compiled with Intel's C++ compiler ($\mathsf{icpc}$ version 11.0).


\begin{figure}[!t]
    \centering
    \subfloat[$L_1$ accuracy]{\includegraphics[width=0.48\hsize]{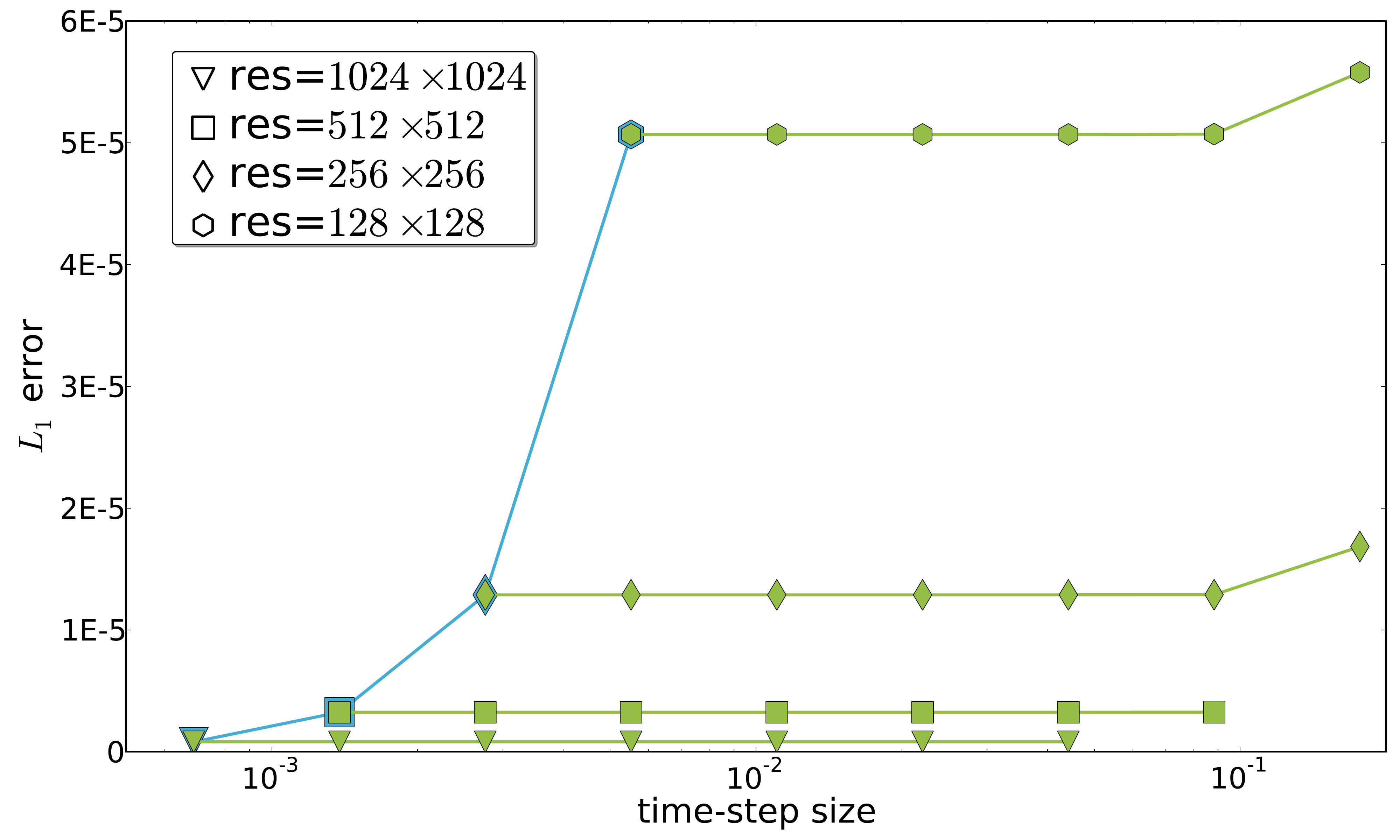}\label{fig:t00_err}}
    \quad
    \subfloat[$L_1$ cost]{\includegraphics[width=0.48\hsize]{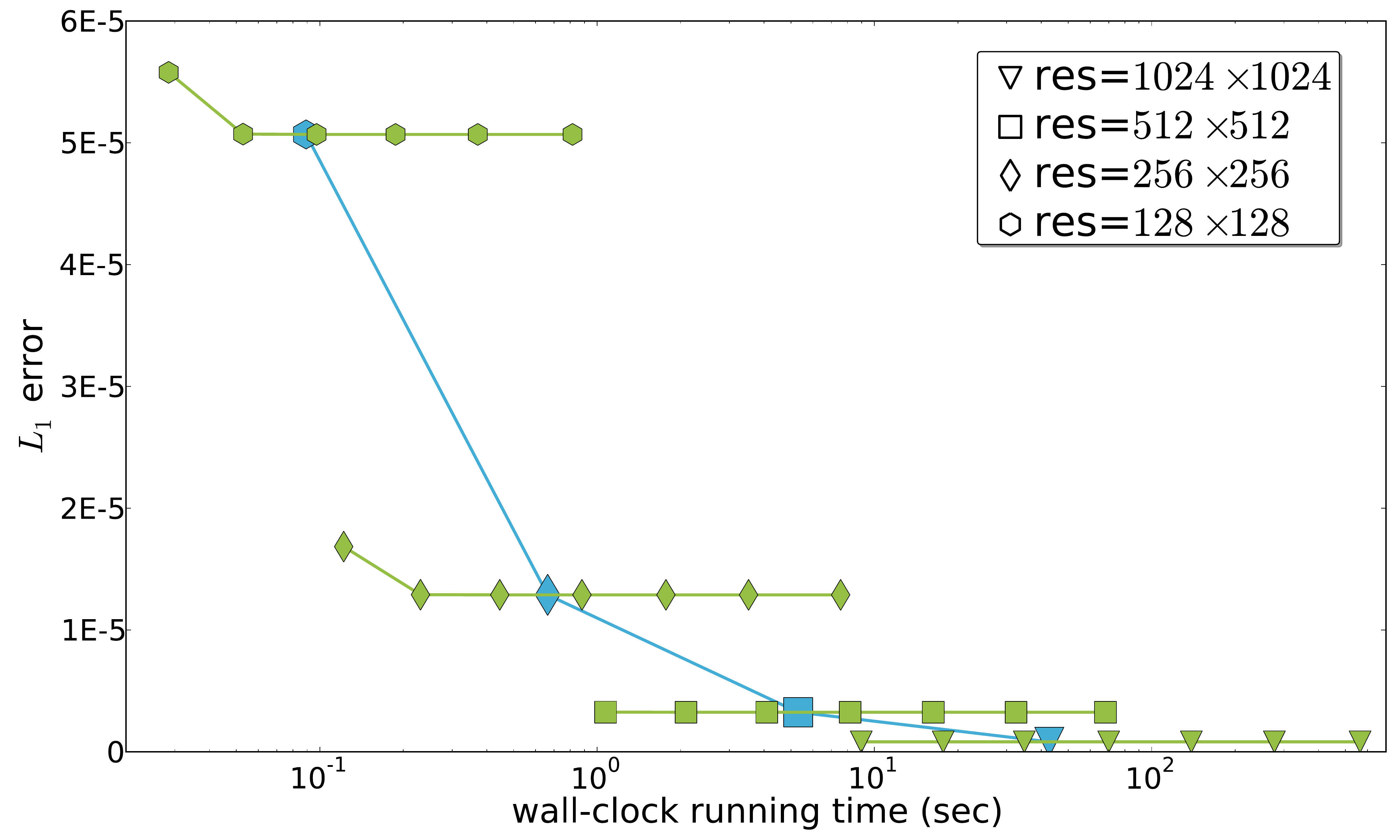}\label{fig:t00_time}}
    \\
    \subfloat[$L_\infty$ accuracy]{\includegraphics[width=0.48\hsize]{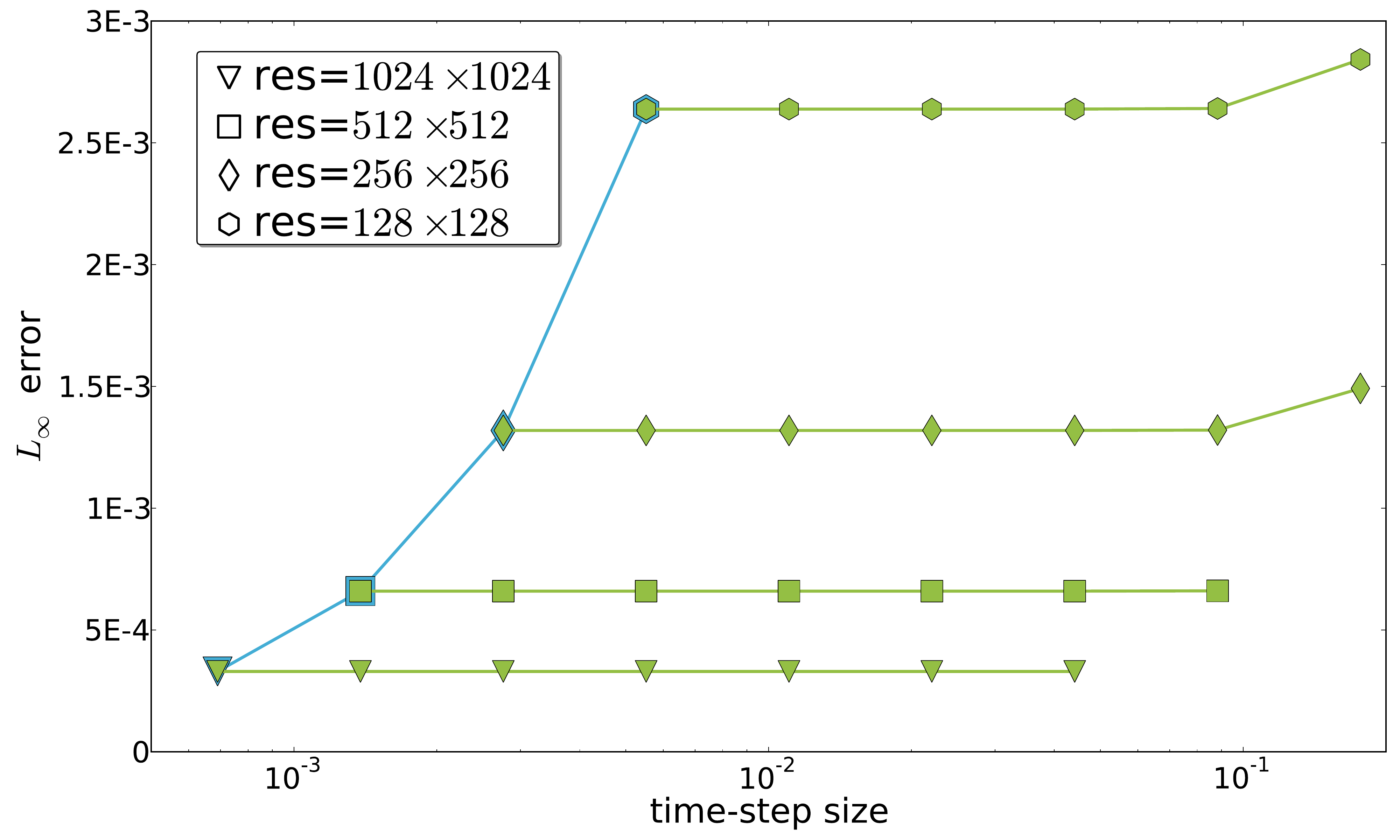}\label{fig:t00_err_inf}}
    \quad
    \subfloat[$L_\infty$ cost]{\includegraphics[width=0.48\hsize]{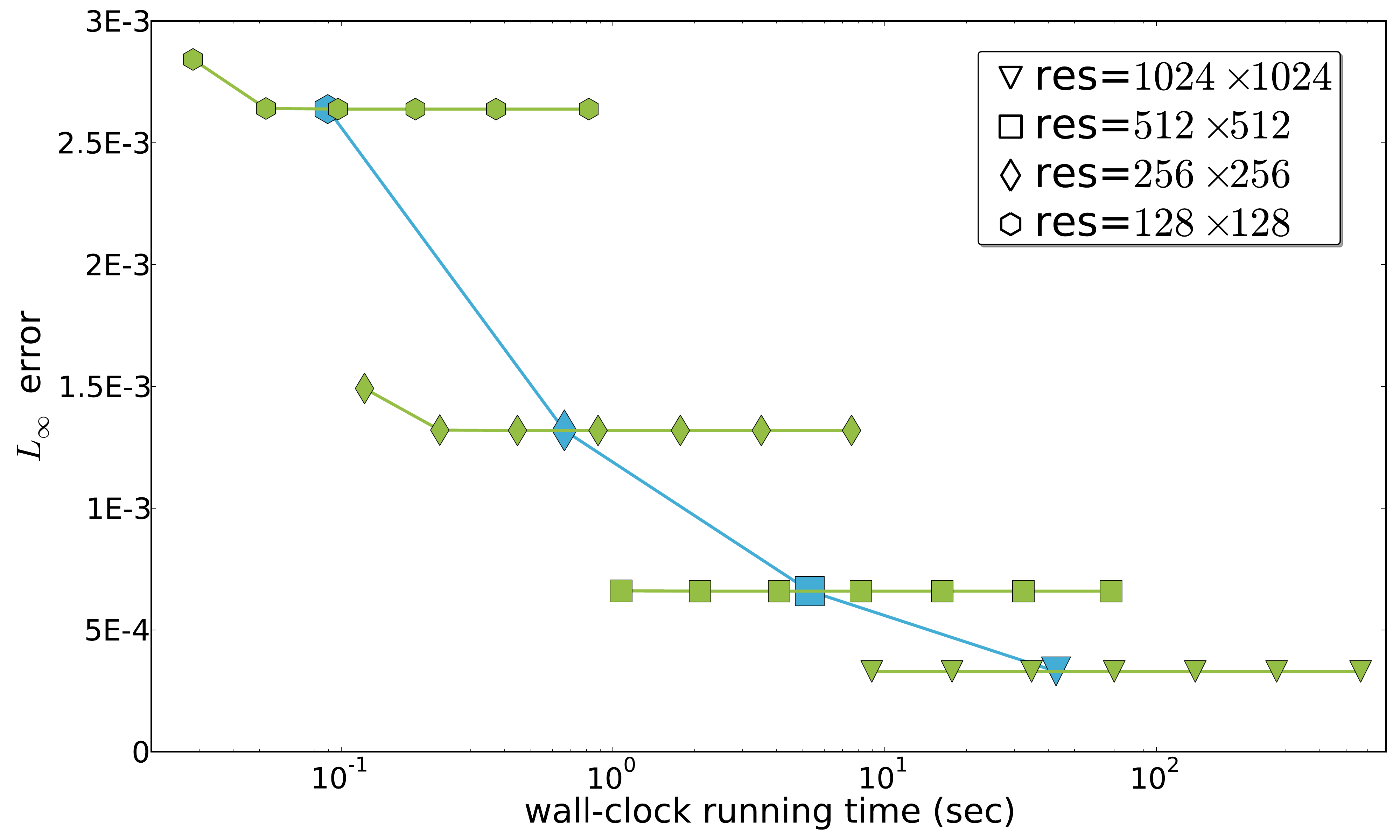}\label{fig:t00_time_inf}}
    \caption{ {\bf Accuracy and computation cost of experiment~1:} the first row uses $L_1$ error metric, and the second
    row uses $L_\infty$ error metric. (a) and (c) plot error curves in terms of different resolution. (b) and (d) show
    error curves in terms of wall-clock computation time.\label{fig:t00}}
\end{figure}
\paragraph{Experiment 1:}
We start with a simple homogeneous min-cost control problem 
with zero boundary/terminal conditions. Namely,
\begin{equation}
\begin{split}
    u_t(\xx,t) - |\nabla u(\xx,t)| = -1,\qquad & \xx\in\Omega, t\in[0,T] \\
    u(\xx,t) = 0, \qquad & \xx\in\partial\Omega,t\in[0,T] \textrm{ or }\xx\in\Omega, t=T
\end{split}
\end{equation}
The analytic solution is $u(\xx,t) = \min\{\phi(\xx), T-t\}$, where $\phi(\xx)$ is the distance
from $\xx$ to $\partial\Omega$, i.e.~$\phi(\xx)=\min\{ x, y, 1-x, 1-y \}$.

In Figure~\ref{fig:t00}, we compare the accuracy and efficiency of the explicit and implicit methods
using $T=1.2$.
The format of Figure~\ref{fig:t00_err} is similar to the linear examples in Figure~\ref{fig:hybrid_vs_imp_vs_exp}.
The blue curve presents explicit method results for various spatial resolutions
(encoded by different marker shapes) with the timestep sizes specified by the CFL.
Each green curve corresponds to the implicit method results for each of those resolutions. We start
from the explicit method's timestep (the leftmost marker) and repeatedly double it as we move to
the right along that green curve.

\begin{remark}\label{rem:cost_update}
    Unlike the linear 1D examples considered in section~\ref{s:advection}, here it is not enough to know
    that some of the green markers fall below the blue accuracy-versus-$h$ curve for a fair comparison
    of these methods' efficiency. There are three reasons for this:
    \begin{itemize}
        \item Since $\domain$ is 2-dimensional, each step to the right along the blue curve decreases
            the number of gridpoint-updates by the factor of 8,
            while each step to the right along any green curve
            yields the factor of 2 reduction only.
            See also Remark~\ref{rem:imp_cost_higher_dim}.
        \item Our modified fast-marching implementation has complexity $O(M\log M)$ per time
            slice, unlike the $O(M)$ per time slice complexity in the explicit approach.
        \item In addition to the computational complexity comparison, explicit update
            equation~\eqref{eq:explicit_discr} is linear and thus less costly than the quadratic
            update used in the implicit discretization.
    \end{itemize}
    Thus, a fair comparison should be based on \emph{``accuracy vs.\ the running time''} analysis;
    see Figure~\ref{fig:t00_time}.
\end{remark}

A similar accuracy/efficiency comparison based on $L_\infty$ errors is provided in Figures~\ref{fig:t00_err_inf} and \ref{fig:t00_time_inf}.
It shows that the presence of shock lines does not impede the advantages of implicit approach.

This simple example might seem counter-intuitive, since the speed $f$ is not stiff at all. However,
it is a very special case, since the solution becomes stationary ($u(\bm{x},t)=\phi(\bm{x})$) for
$t<T-0.5$. Thus, both methods can be interpreted as iterative schemes converging to the stationary
solution, which explains the flatness of green accuracy curves. Of course, a steady-state solution
could be obtained more efficiently by solving the stationary PDE directly. 
The next example shows that the implicit methods also have similar performance advantages even for
non-stiff problems, provided the dependence of $u$ on $t$ is nearly linear; see also
Remark~\ref{rem:stiffness_and_more}.







\paragraph{Experiment 2:} We consider a simple test case
\begin{figure}[!t]
    \centering
    \subfloat[accuracy ($\lambda=0.1$)]{\includegraphics[width=0.48\hsize]{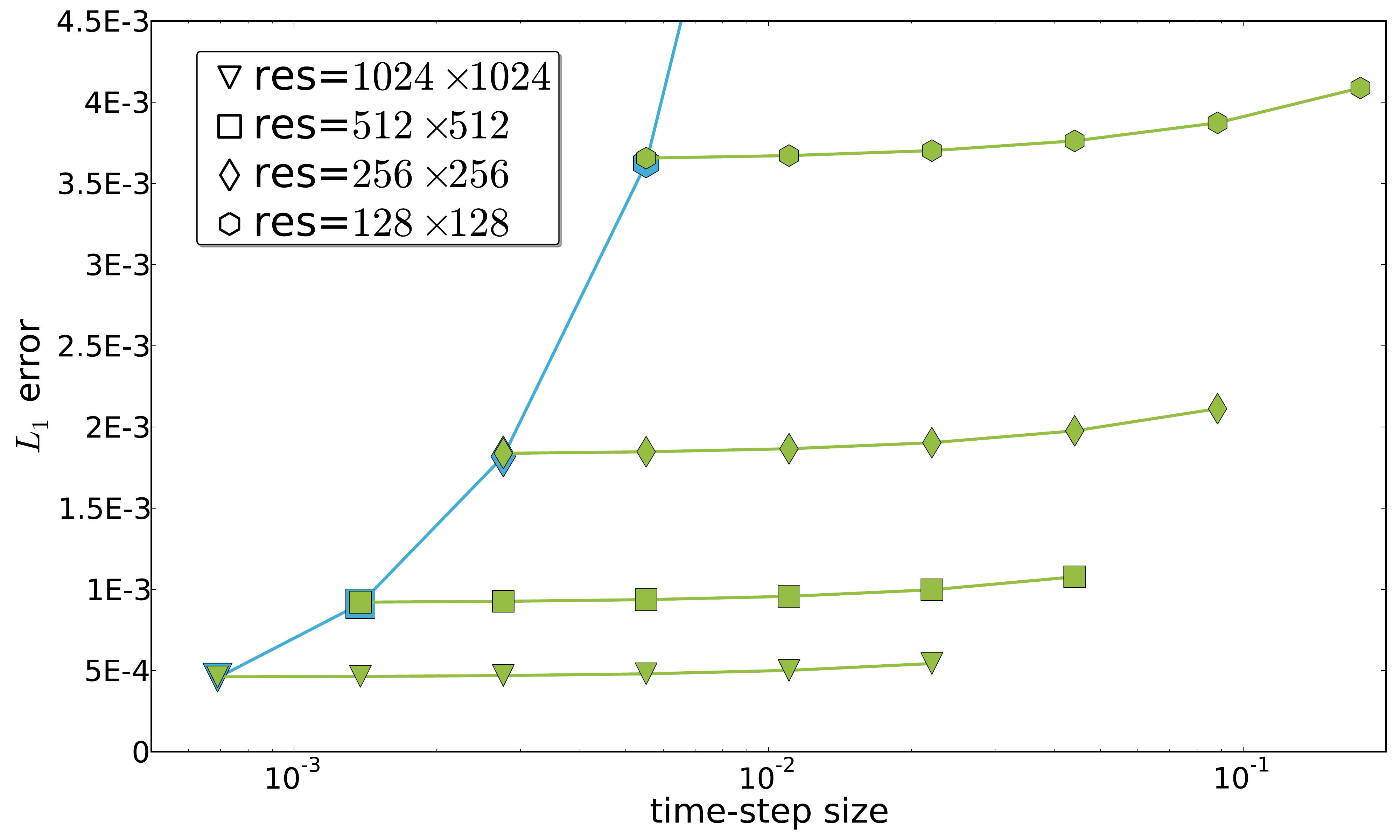}\label{fig:test22_err_0.1}}
    \quad
    \subfloat[cost ($\lambda=0.1$)]{\includegraphics[width=0.48\hsize]{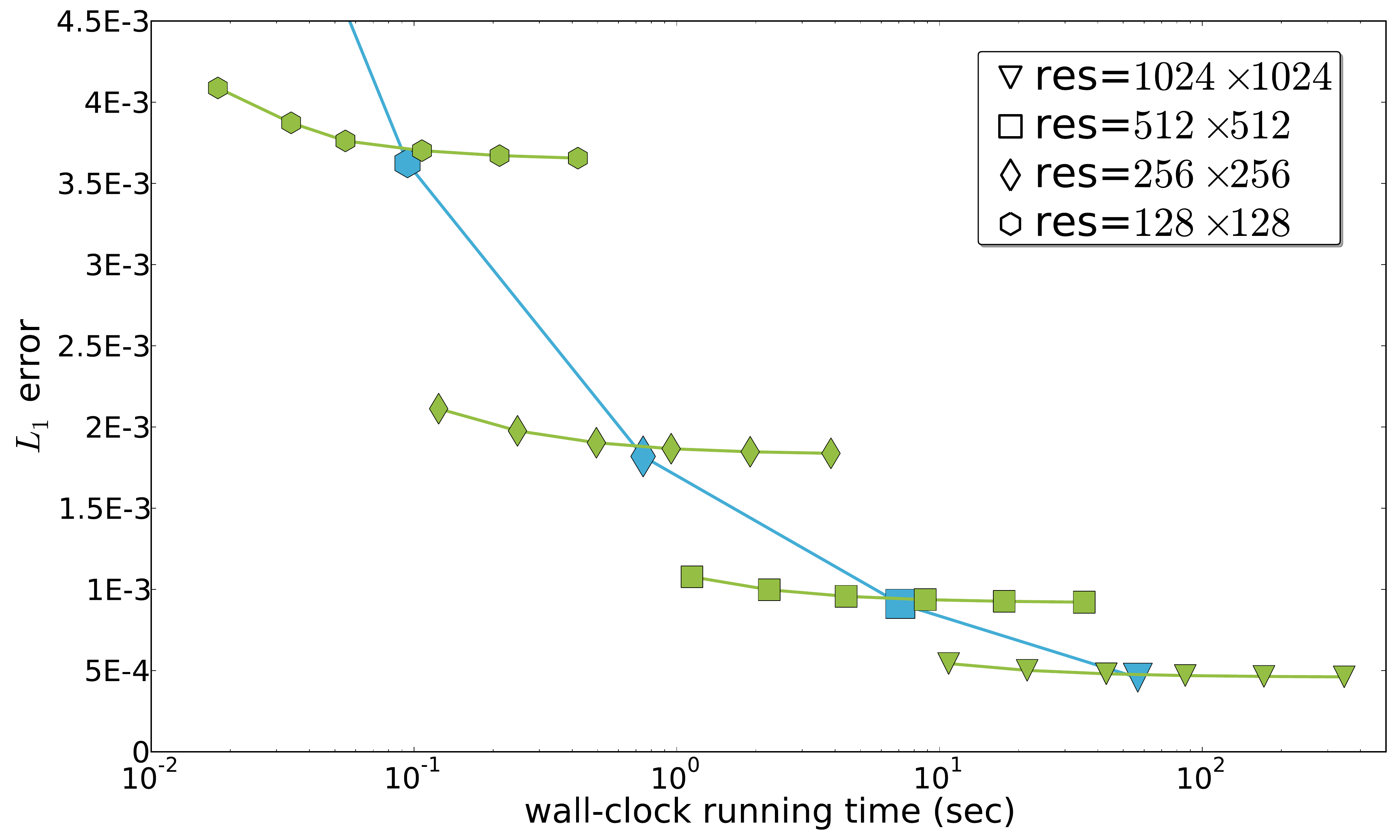}\label{fig:test22_cost_0.1}}
    \\
    \subfloat[accuracy ($\lambda=0.25$)]{\includegraphics[width=0.48\hsize]{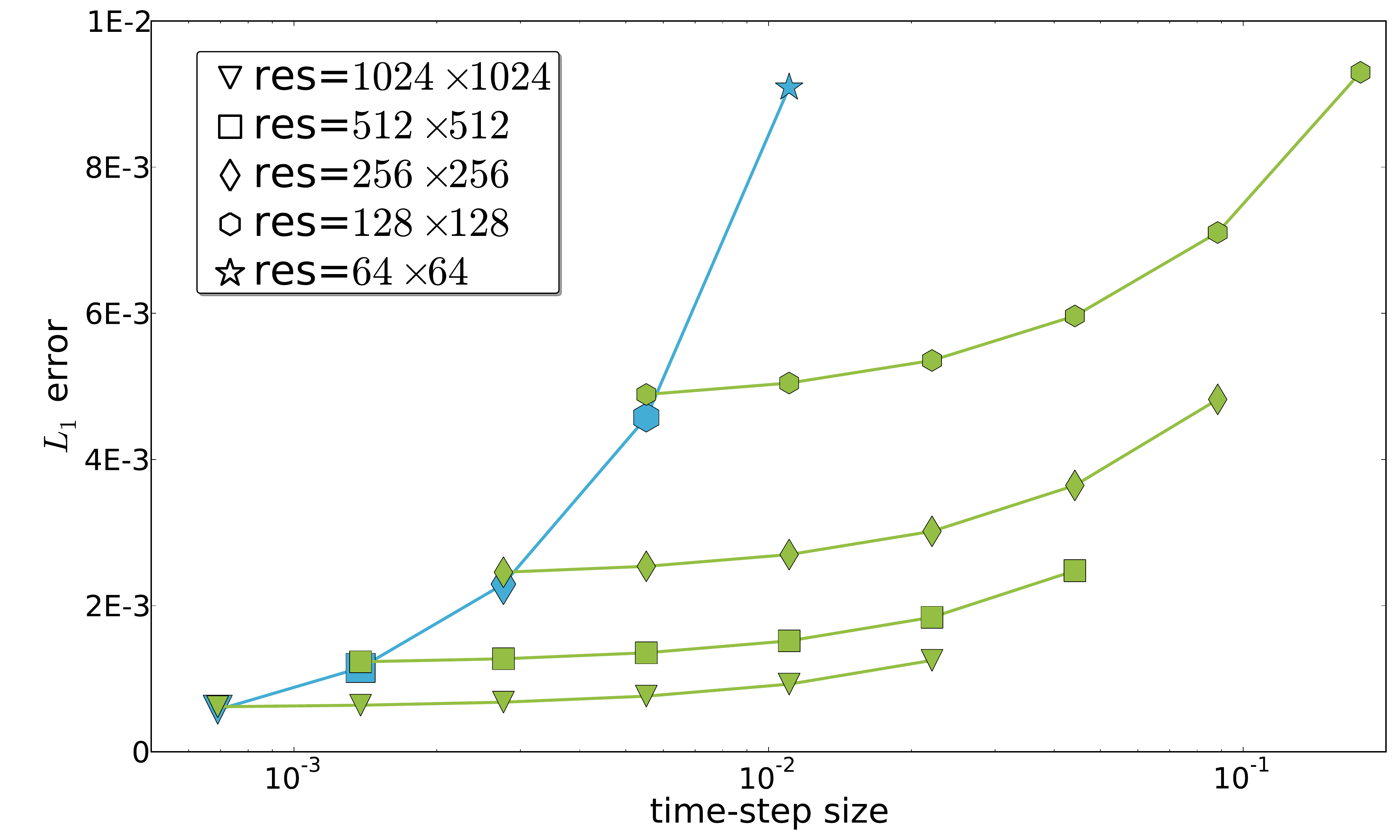}\label{fig:test22_err_0.25}}
    \quad
    \subfloat[cost ($\lambda=0.25$)]{\includegraphics[width=0.48\hsize]{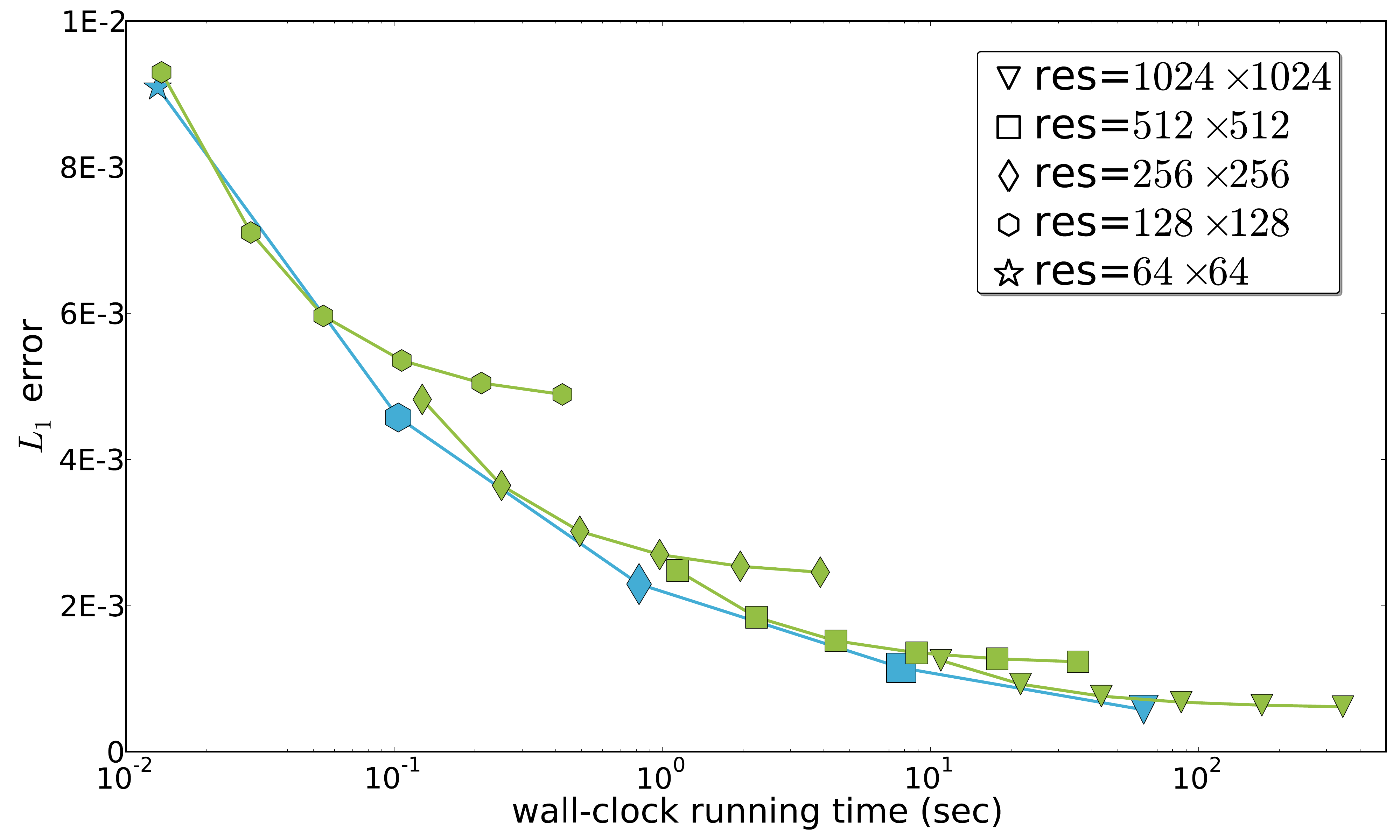}\label{fig:test22_cost_0.25}}
    \\
    \subfloat[accuracy ($\lambda=0.8$)]{\includegraphics[width=0.48\hsize]{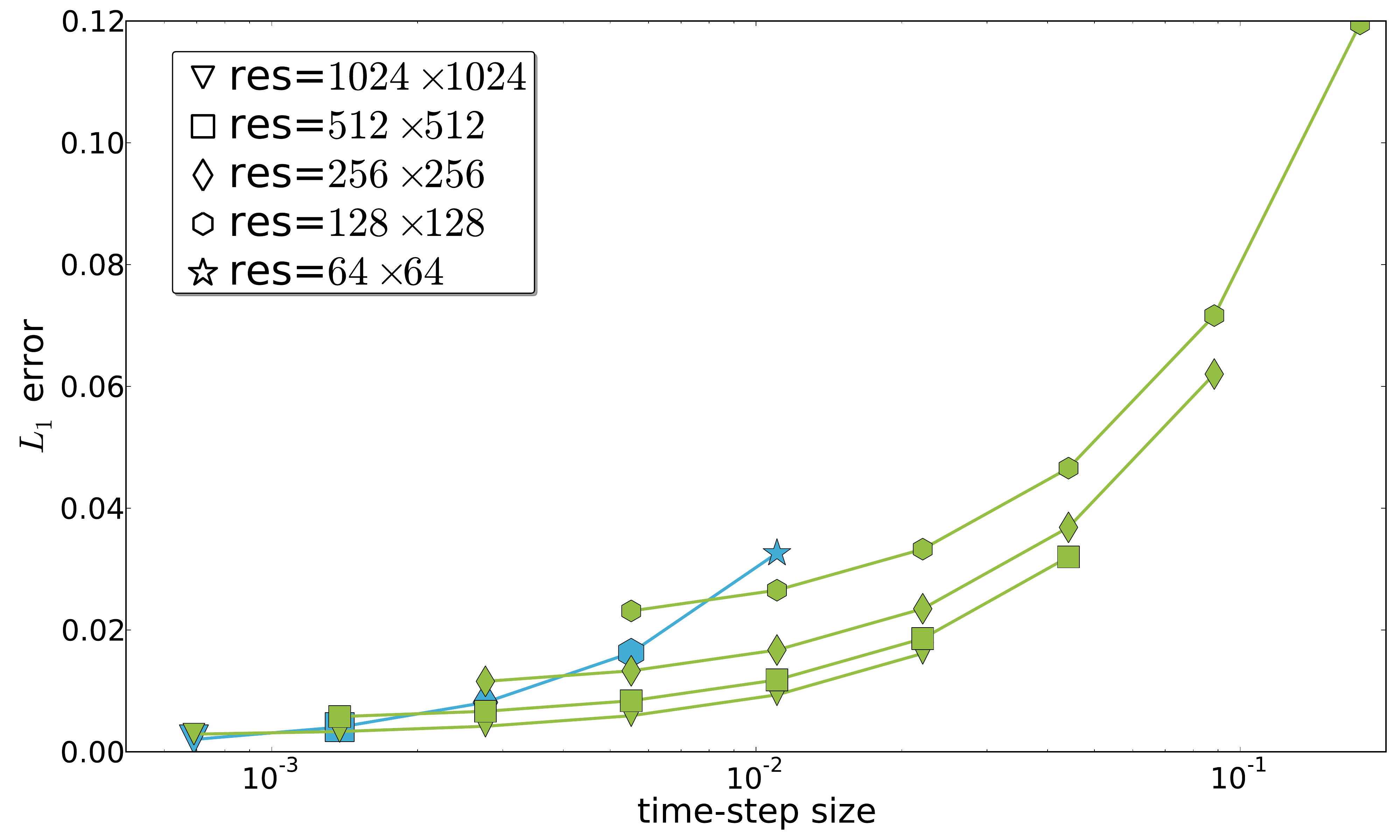}\label{fig:test22_err_0.8}}
    \quad
    \subfloat[cost ($\lambda=0.8$)]{\includegraphics[width=0.48\hsize]{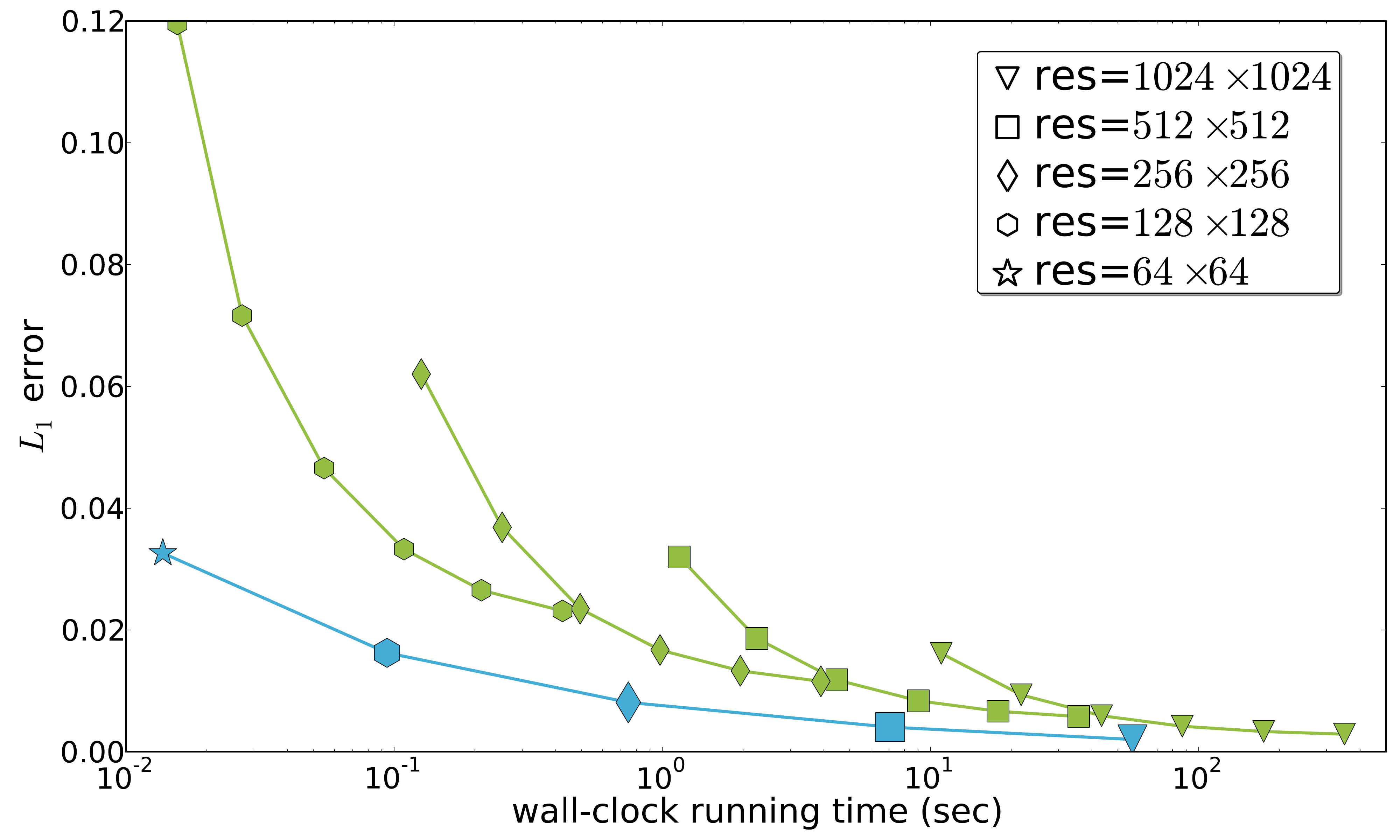}\label{fig:test22_cost_0.8}}
    \caption{ {\bf Accuracy and computation cost of experiment 2:} figures in the first row correspond to the equation with $\lambda=0.1$,
    the second row corresponds to the case of $\lambda=0.25$, and the third rows corresponds to the case of $\lambda=0.8$.
    Figures in the first column plot error curves in terms of different resolution,
    and figures in the second column show error curves in terms of wall-clock computation
    time.\label{fig:test22}}

\end{figure}
\begin{equation}
\begin{split}
    u_t(\xx,t) - \frac{|\nabla u(\xx,t)|}{2y+1} = -1,\qquad & \xx\in\Omega, t\in[0,T] \\
    u(\xx,t) = e^{\lambda t},\qquad & \xx\in[0,1]\times\{0\},\;t\in[0,T] \\
    u(\xx,t) = y+y^2 + e^{\lambda(t+y+y^2)},\qquad& \xx\in\Omega,t=T,
\end{split}
\end{equation}
with $T=1.2$ and the rest of the boundary considered to be ``outflow''.
In this example, all the characteristics are perpendicular to x-axis; therefore, this is
essentially a 1D problem (i.e. $u_x=0$). Its analytic solution is known as
$u(\xx,t) = y + y^2 + e^{\lambda(t+y+y^2)}$. The time-dependence of $u$ is completely due to the
time-dependent boundary conditions, and its effect can be tuned by changing the value of the
parameter $\lambda>0$.  The speed function $f$ is 
only moderately stiff
(its values range from 1/3 to 1).

Figure~\ref{fig:test22} shows that, despite this lack of stiffness, the
implicit method is better for smaller $\lambda$s (when the
dependence of $u$ on $t$ is not too far from linear on $[0,T]$), but clearly loses to the explicit when
$\lambda>0.25$.   This is to be expected, since larger time-steps will fail to capture the rapid changes
in boundary conditions.



\paragraph{Experiment 3:} 
\begin{figure}[!t]
    \centering
    \subfloat[speed function]{\includegraphics[height=32mm]{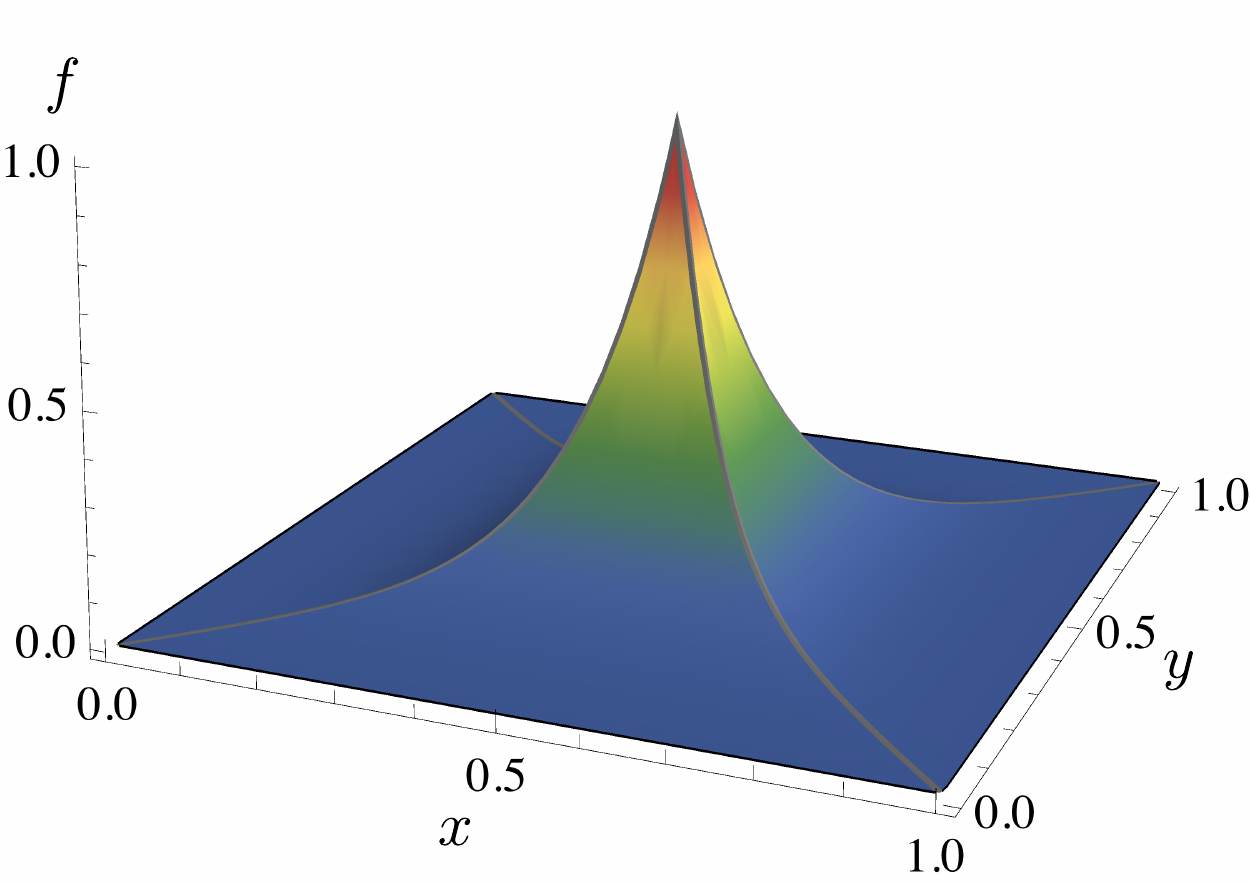}\label{fig:t40_speed}}
    \;
    \subfloat[boundary condition]{\includegraphics[height=32mm]{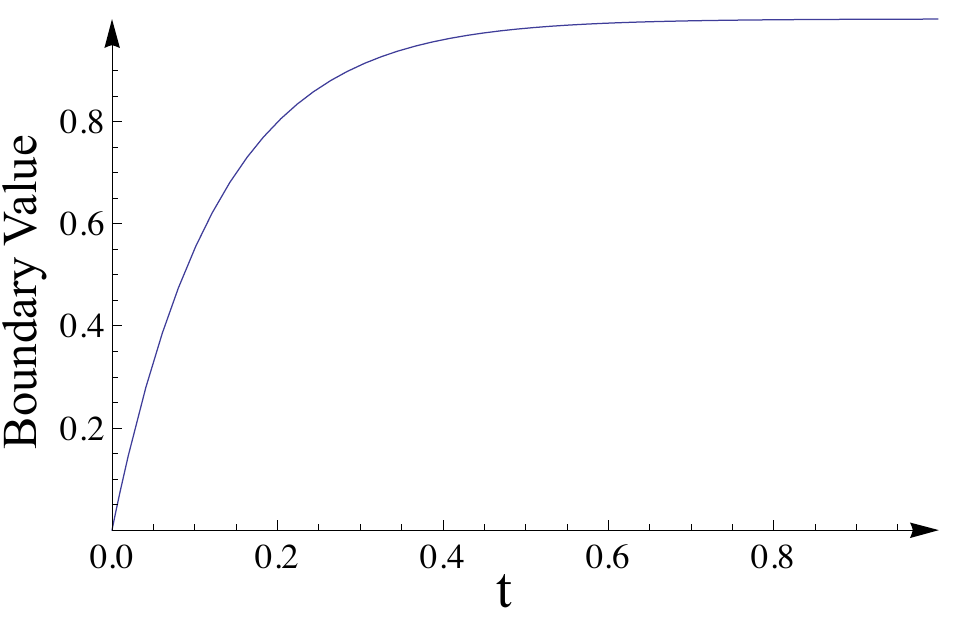}\label{fig:t40_bc}}
    \;
    \subfloat[solution at $t=0$]{\includegraphics[height=32mm]{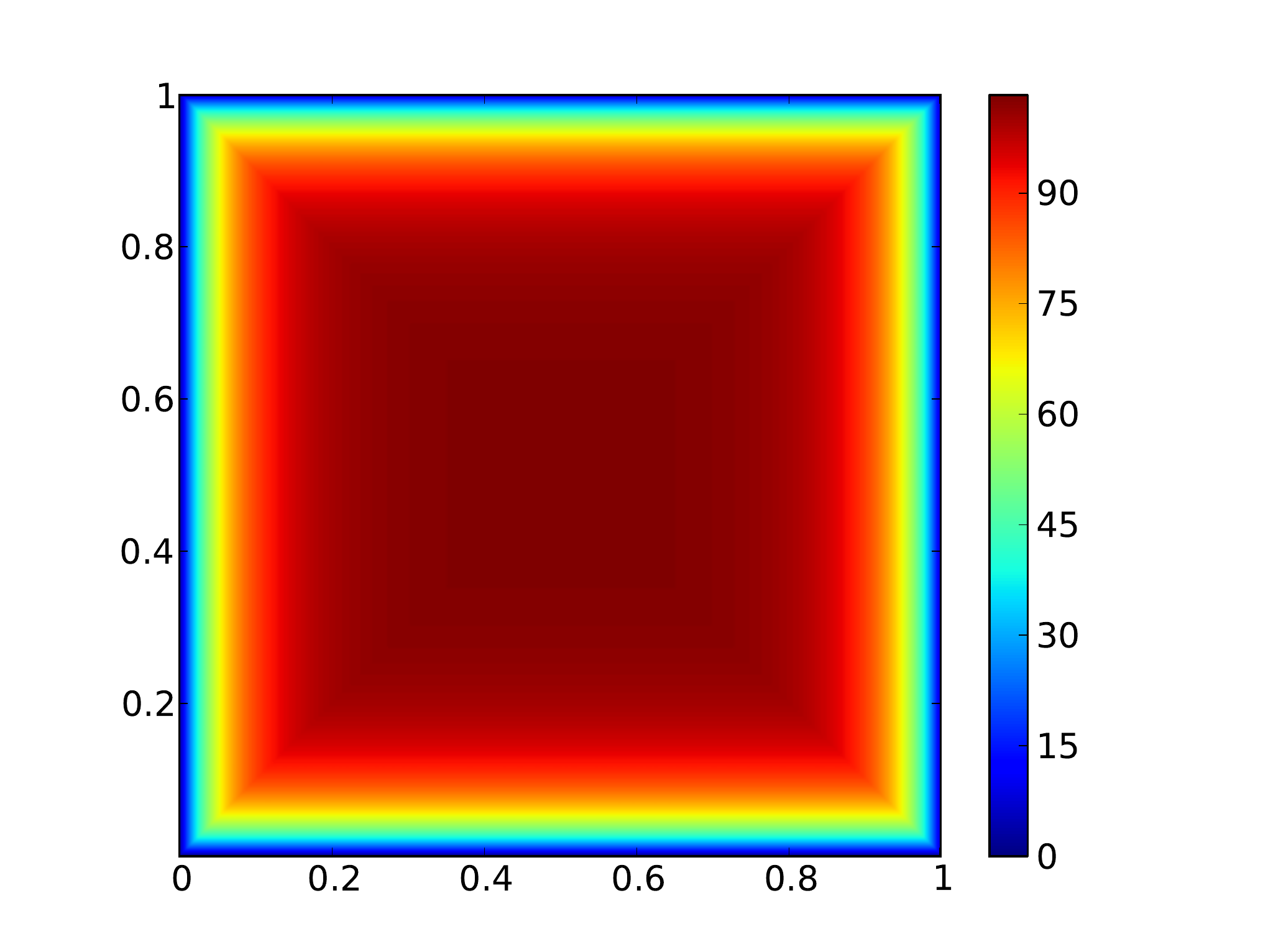}\label{fig:t40_t0}}
    \caption{ (a) A stiff speed function for which large speed values only appear in a small region at the center of the domain.
    The rest of the domain has small speed value close to $\frac{1}{2^\gamma}$. (b) The profile of spatially invariant boundary condition.
    (c) A color-mapped plot of the solution $u(\xx,0)$ with $\gamma=11$.
    }
\end{figure}
\begin{figure}[!t]
    \centering
    \subfloat[accuracy ($\gamma=11$)]{\includegraphics[width=0.48\hsize]{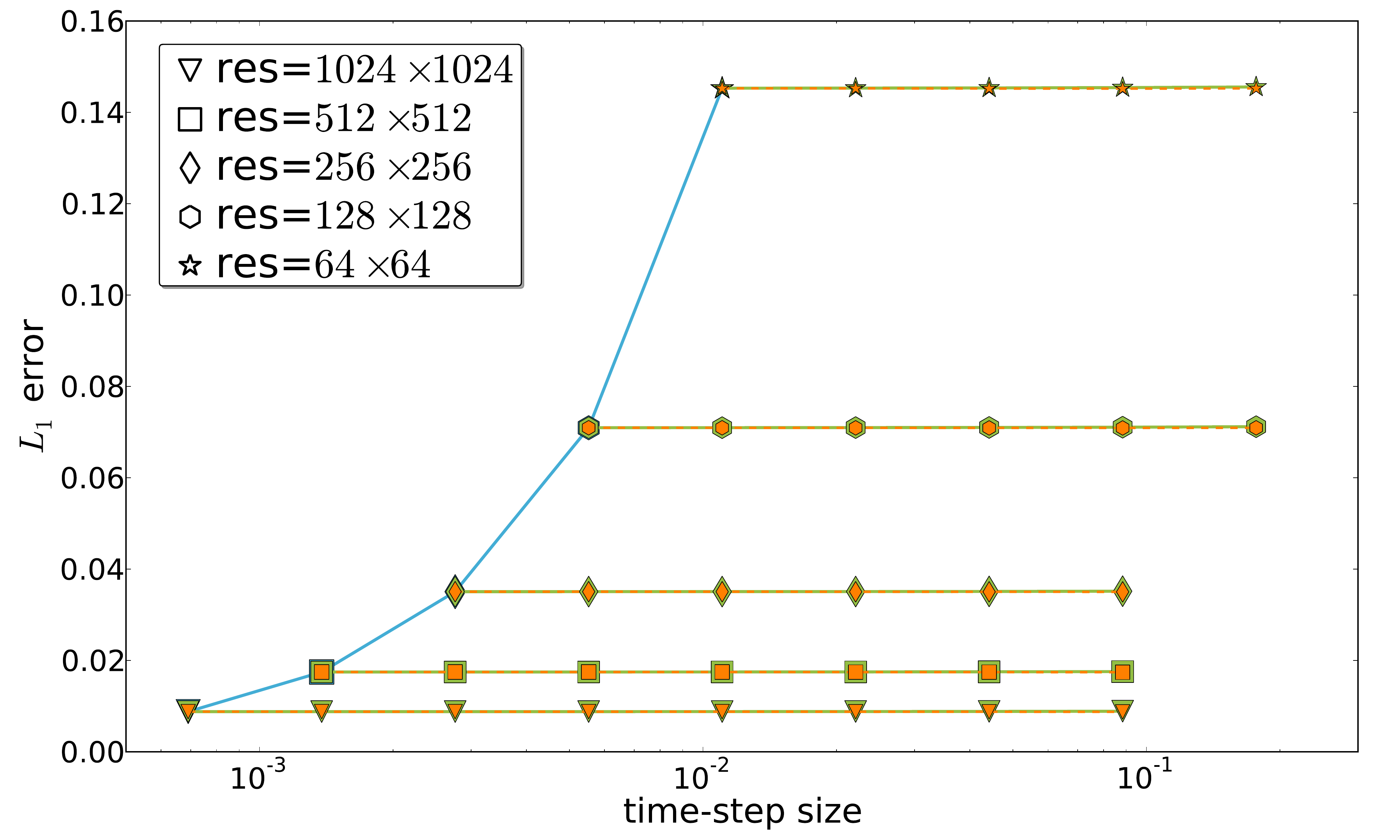}\label{fig:test40_err}}
    \quad
    \subfloat[cost ($\gamma=11$)]{\includegraphics[width=0.48\hsize]{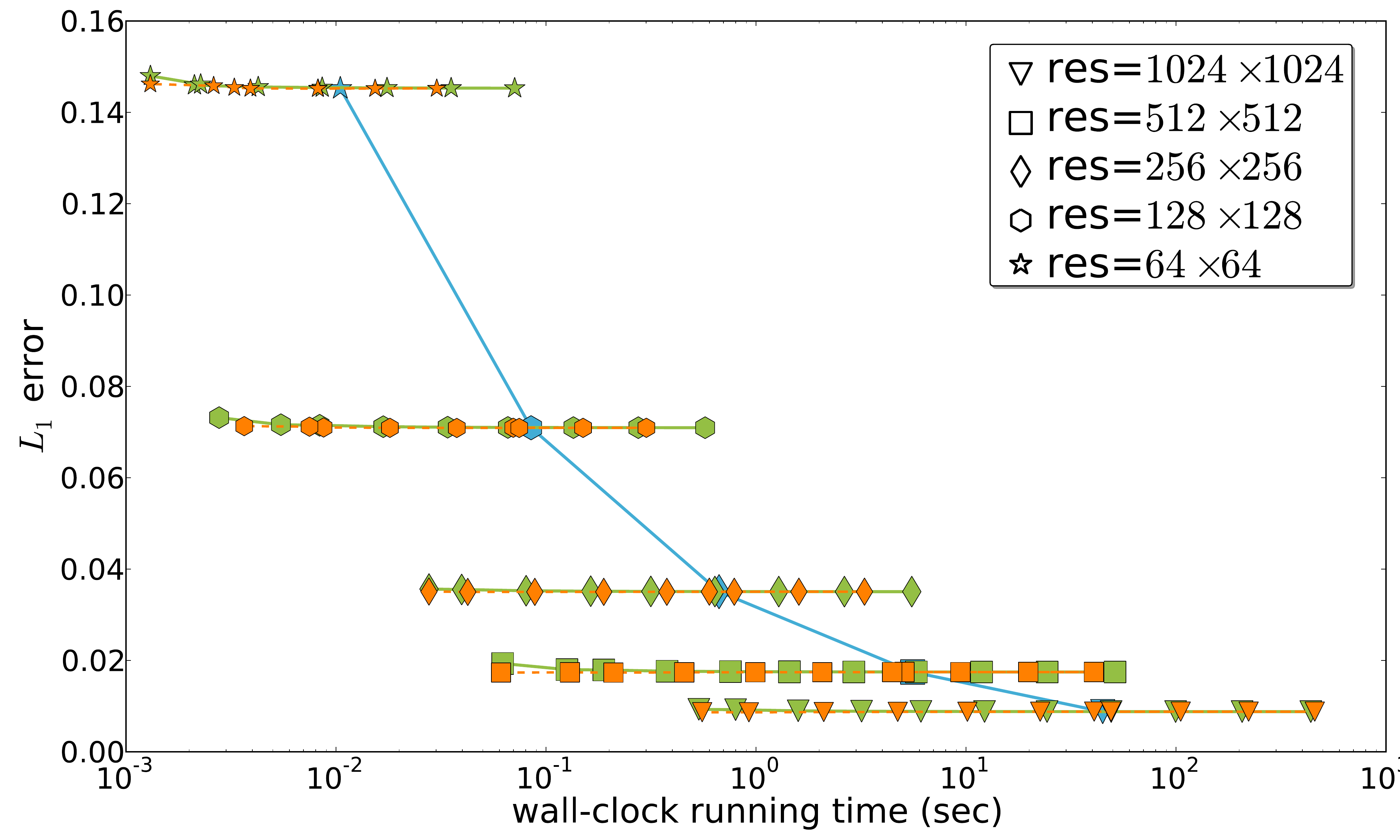}\label{fig:test40_cost}}
    \\
    \subfloat[accuracy ($\gamma=5$)]{\includegraphics[width=0.48\hsize]{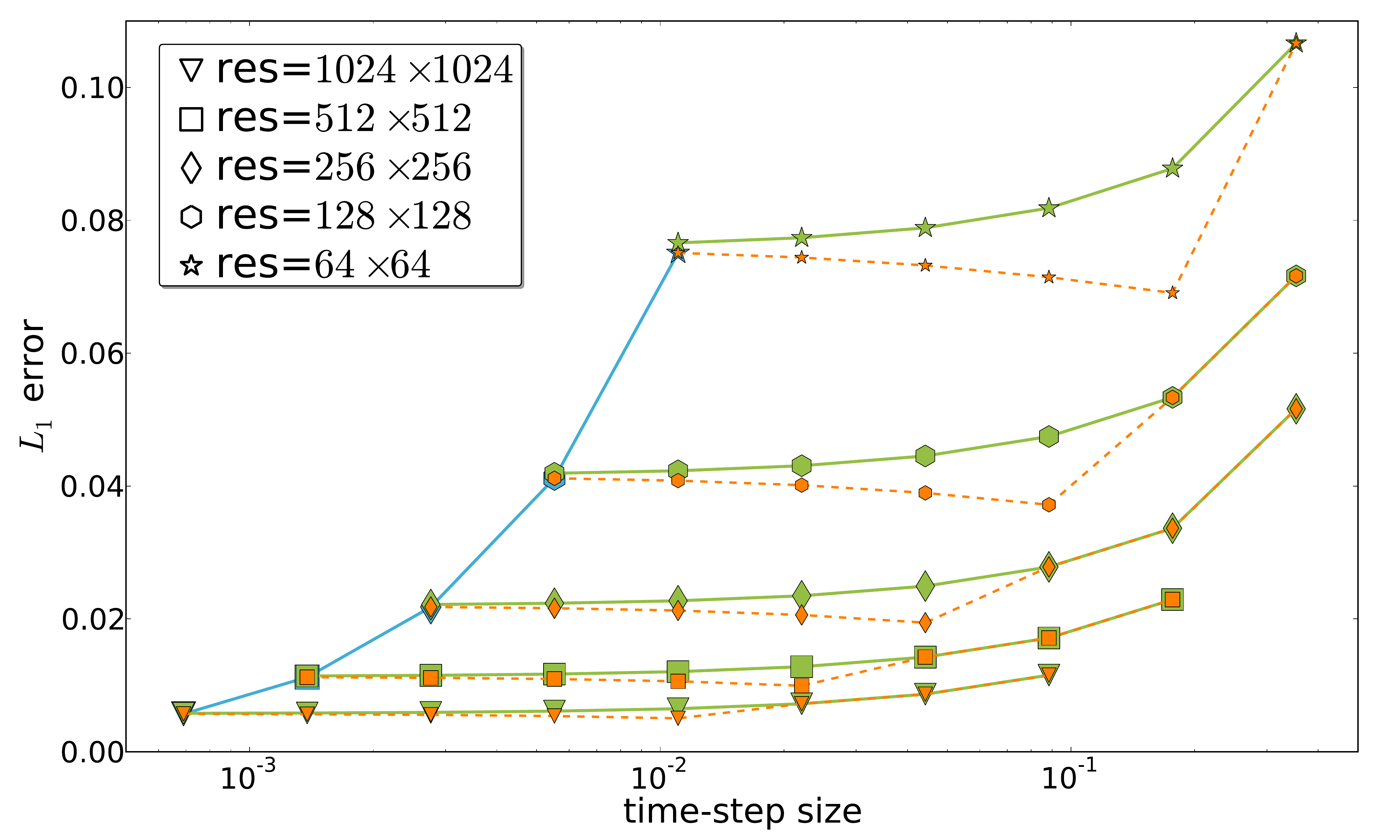}\label{fig:test40_b5_err}}
    \quad
    \subfloat[cost ($\gamma=5$)]{\includegraphics[width=0.48\hsize]{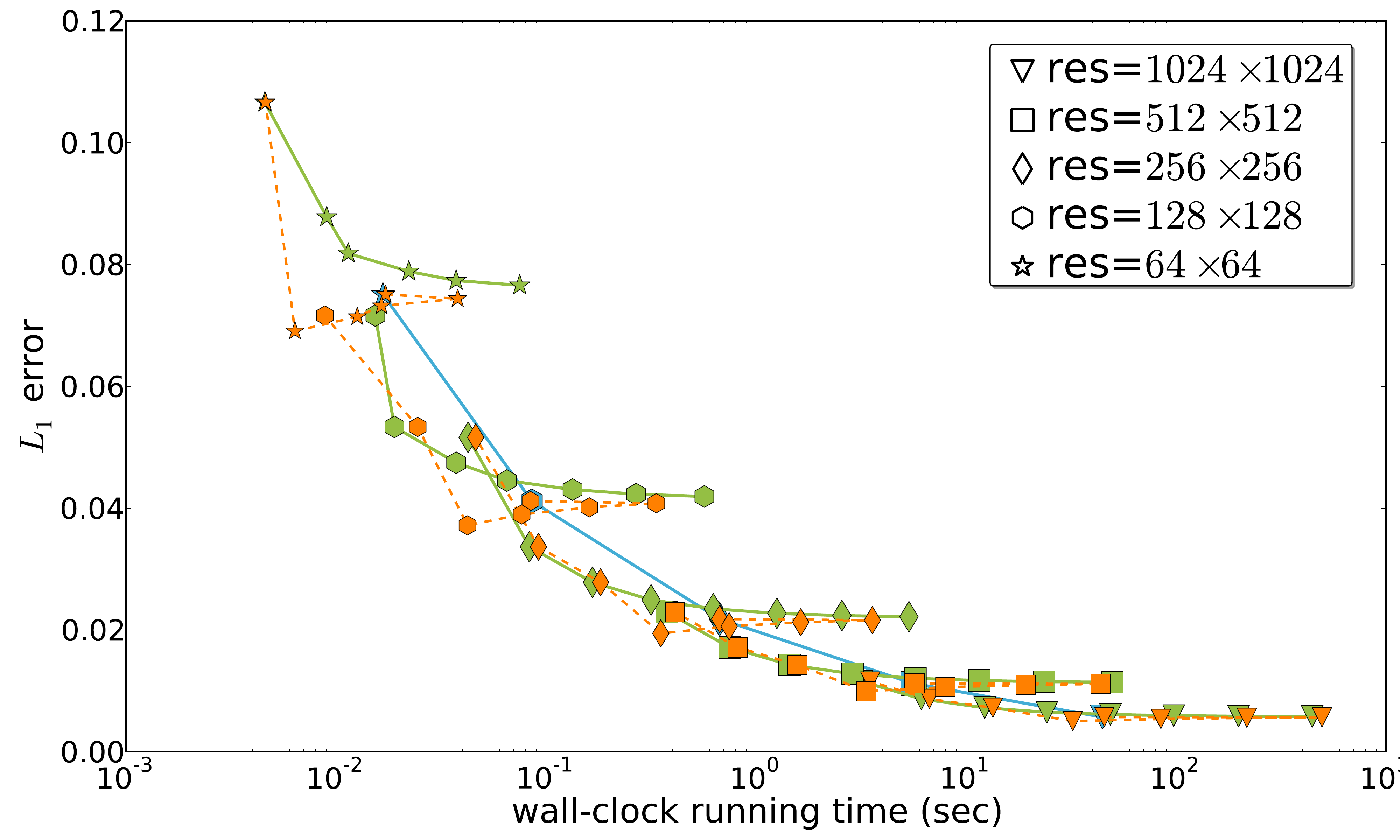}\label{fig:test40_b5_cost}}
    \caption{ {\bf Accuracy and computation cost of experiment 3:} the first row shows the plots of
    error-resolution and error-cost curves with $\gamma=11$;
    the second row shows the same sets of curves with $\gamma=5$.\label{fig:test40}}
\end{figure}
We now turn to stiff examples, and consider a test case
\begin{equation}
    u_t(\xx,t) - \left(\frac{1+2\phi(\xx)}{2}\right)^{\gamma}|\nabla u(\xx,t)| = -1,\qquad \xx\in\Omega, t\in[0,T]
\end{equation}
where $T=1$, $\phi(\xx)$ is the distance from $\xx\in\Omega$ to $\partial\Omega$, 
and $\gamma$ is a parameter to control the stiffness of the speed function ranging in $[0,1]$.
The boundary condition on $\partial\Omega \times [0,T]$ (see Figure~\ref{fig:t40_bc}) is
spatially invariant, defined as
$$ q(\xx, t) = \frac{e^8 - e^{8(1-t)}}{e^8-1}. $$
From the perspective of optimal control, the solution $u(\xx,t)$ is the earliest time of exiting
$\Omega$ if we start from $\xx$ at the time $t$ and pay exit time-penalty $q$ on $\partial\Omega$.
It is easy to show that 
the minimal time to the boundary can be explicitly written as
$$ \tau(\xx) = \frac{2^{\gamma-1}}{\gamma-1}\left(1 - \frac{1}{(1+2\phi(\xx))^{\gamma-1}}\right).$$
Therefore, the analytic solution is
$$ u(\xx,t) = \tau(\xx) + q(t + \tau(\xx)), $$
which we also use to specify the terminal condition
(i.e., $u(\xx,T)=\tau(\xx)+q(T+\tau(\xx))$).

Figures~\ref{fig:test40_err} and \ref{fig:test40_cost} show that the implicit method is much better
than the explicit when we use a stiff speed function with $\gamma=11$. This is due to the fact
that large speeds occur at a small central part of the domain only; thus, the large timesteps can be
used without much impact on accuracy.
Here we also show
the performance of the hybrid method in orange curves, but it yields no advantage over the implicit
in this example. On the other hand, for a lower level of stiffness ($\gamma=5$,
Figures~\ref{fig:test40_b5_err} and \ref{fig:test40_b5_cost}), the advantage of implicit method is
less significant, but the hybrid method offers additional boost in both accuracy and efficiency.
A careful look at orange curves in Figure~\ref{fig:test40_b5_cost} shows complex/non-monotone
behavior of hybrid methods. A similar phenomenon is explained in detail in the next example.




\paragraph{Experiment 4:}
\begin{figure}[!t]
    \centering
    \subfloat[$u(\xx,0)$]{\includegraphics[height=40mm]{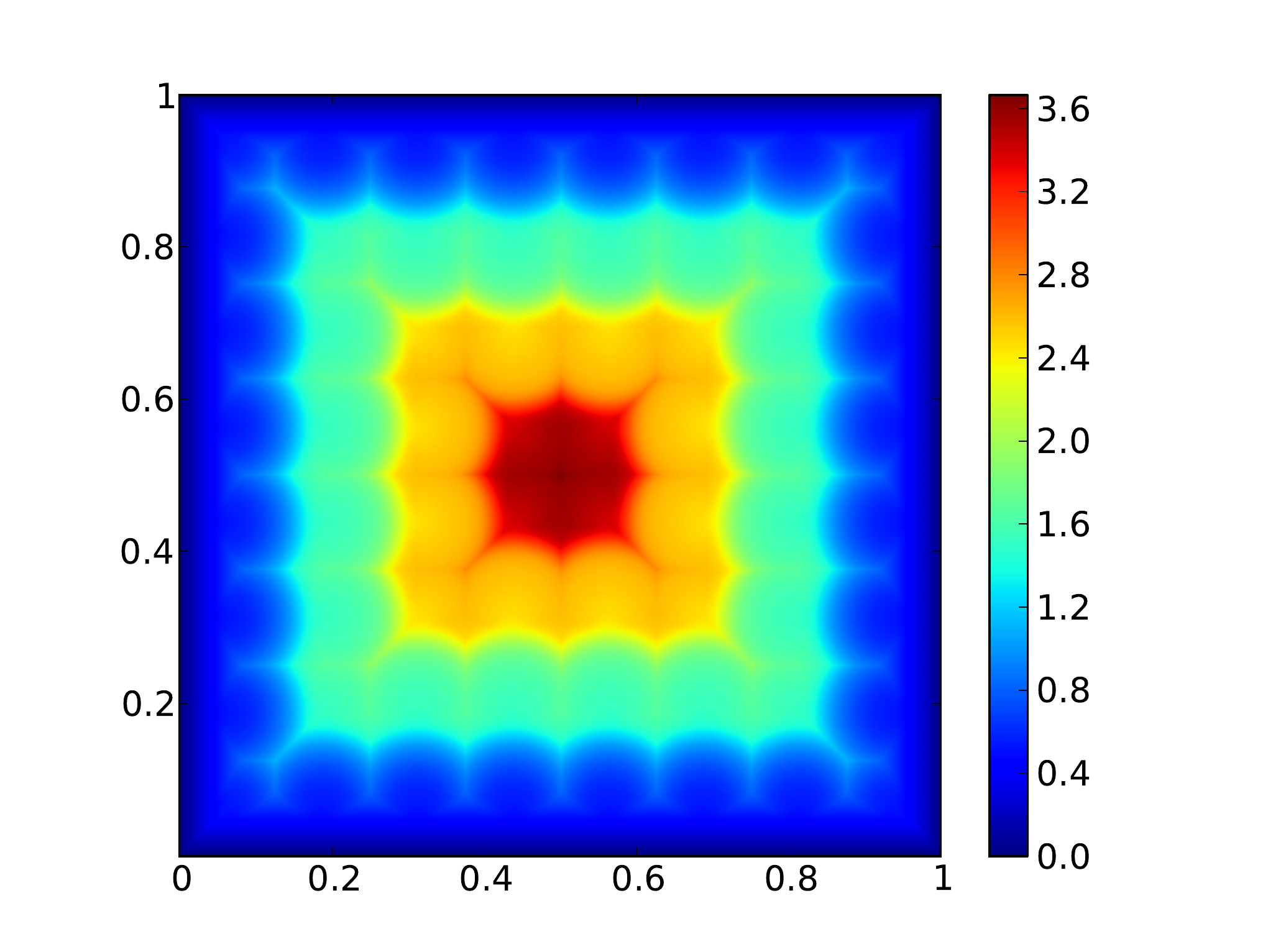}\label{fig:t07_0}}
    \hspace{-8mm}
    \subfloat[$u(\xx,0.2)$]{\includegraphics[height=40mm]{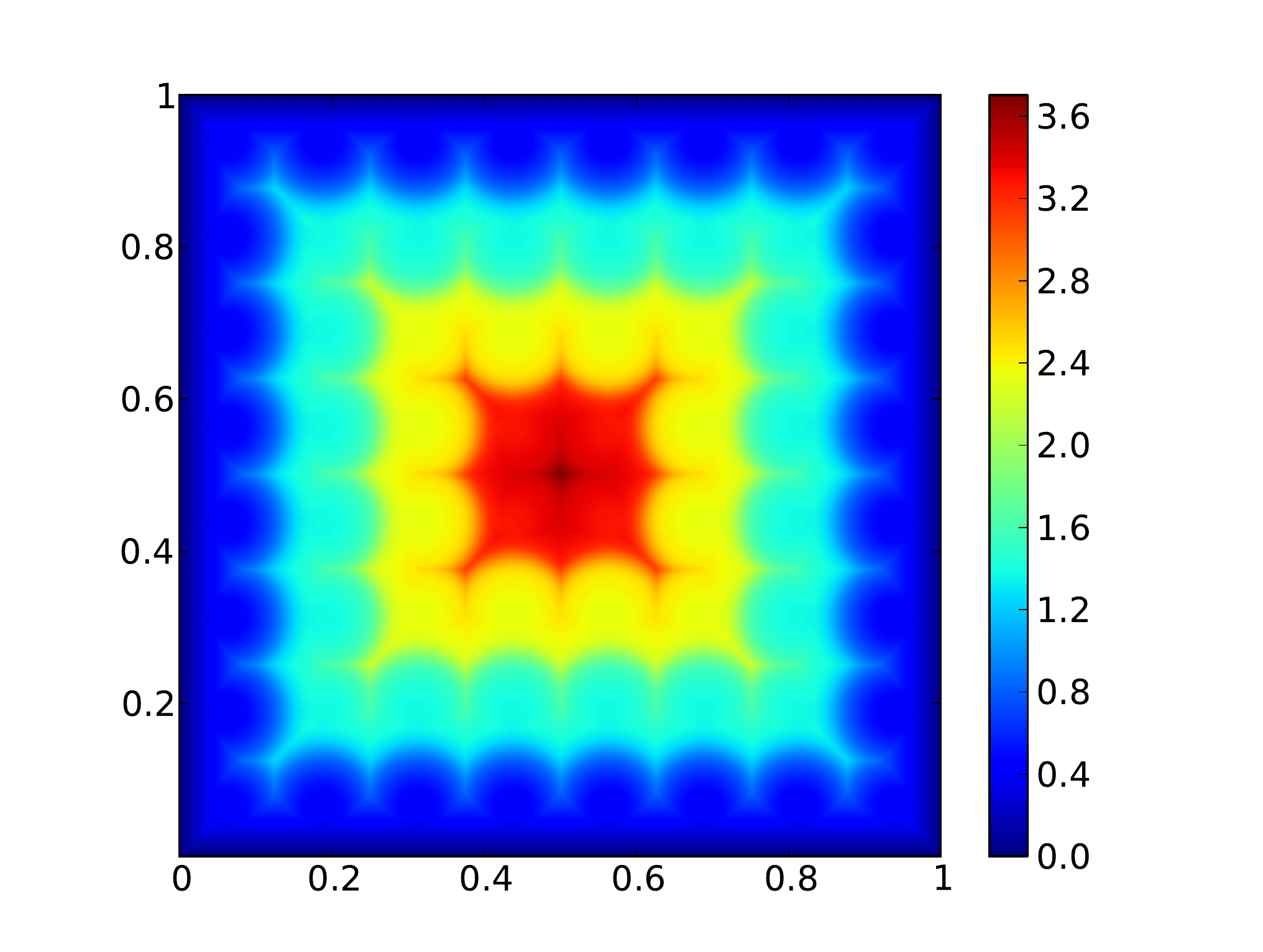}\label{fig:t07_0.2}}
    \hspace{-8mm}
    \subfloat[$u(\xx,0.4)$]{\includegraphics[height=40mm]{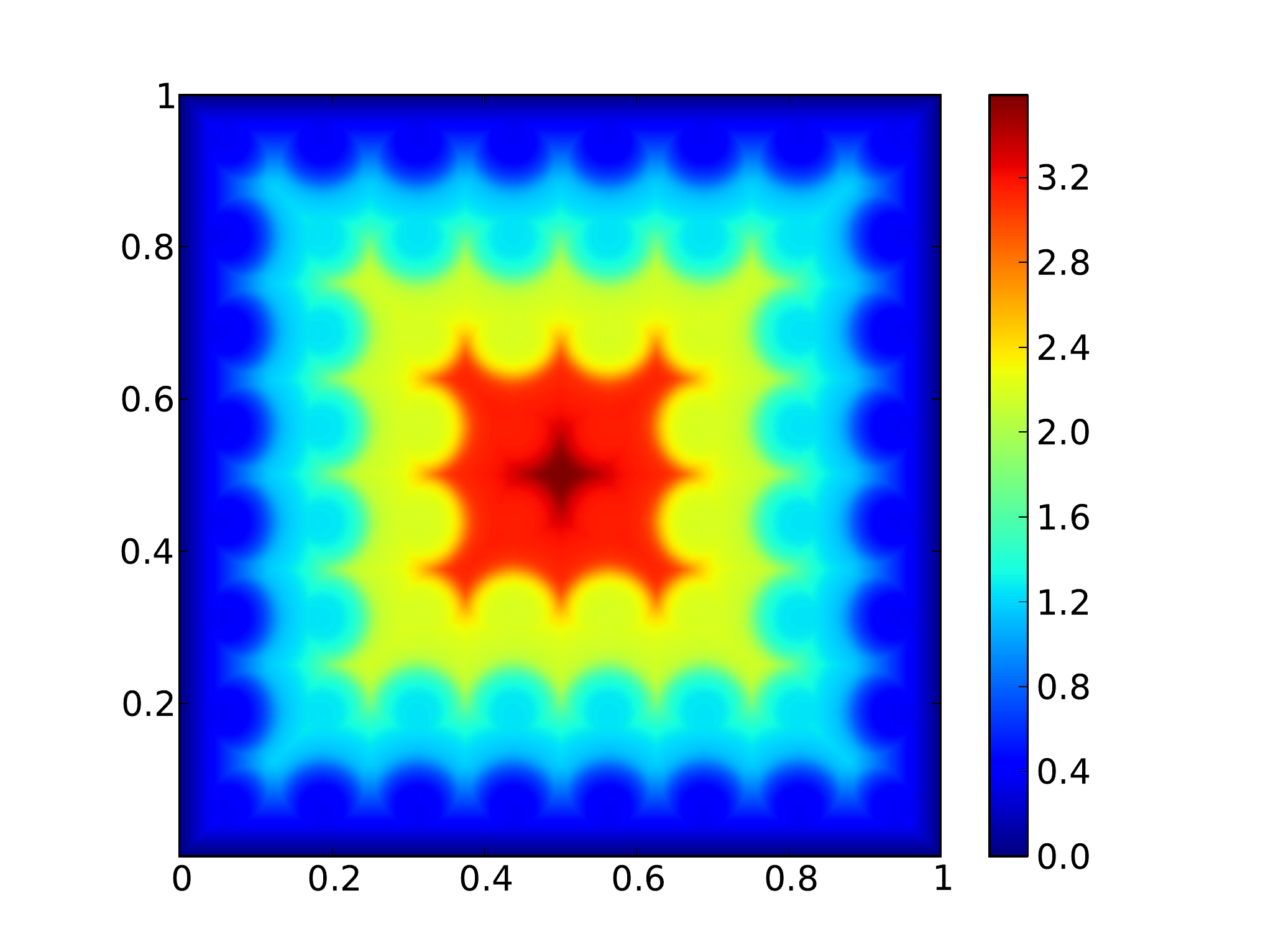}\label{fig:t07_0.4}}
    \caption{ {\bf Color-mapped plot of $u(\xx,t)$ of experiment 4}\label{fig:test07_ret}}
\end{figure}
\begin{figure}[!t]
    \centering
    \subfloat[$L_1$ accuracy]{\includegraphics[height=31mm]{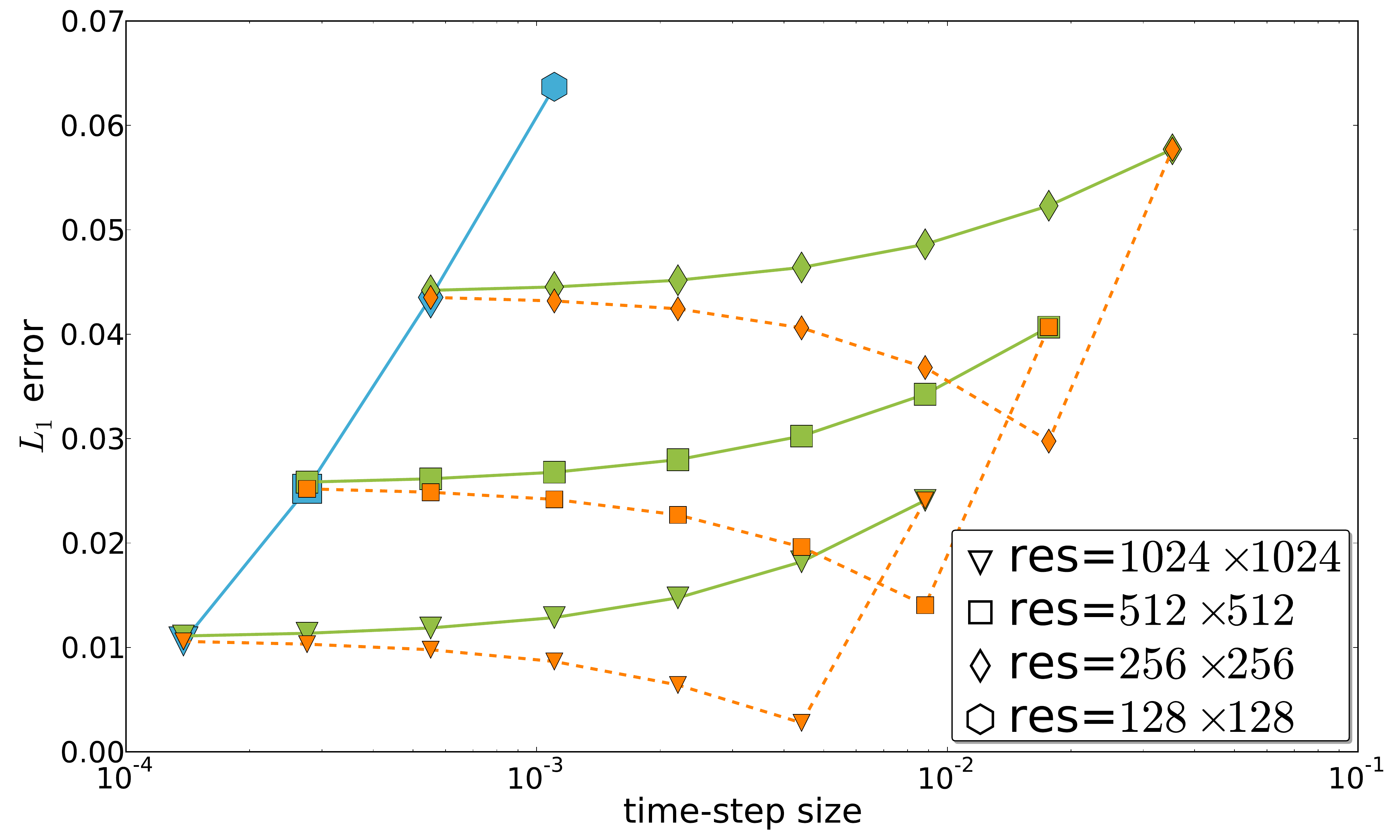}\label{fig:test07_err}}
    \;
    \subfloat[$L_1$ cost]{\includegraphics[height=31mm]{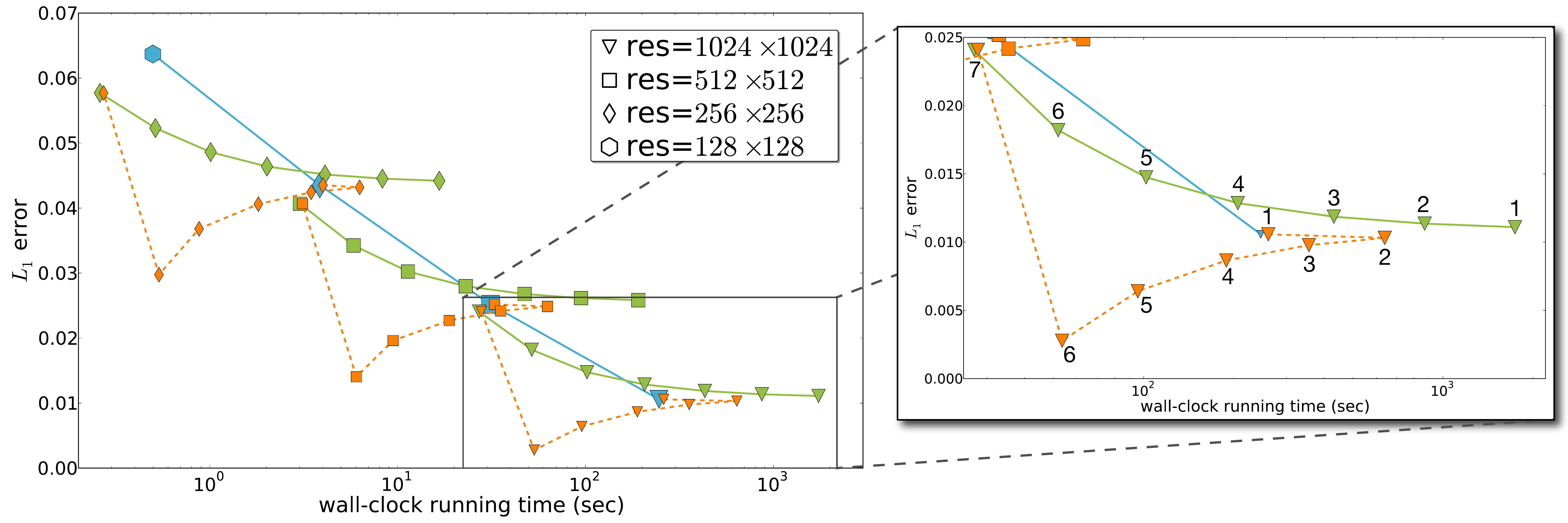}\label{fig:test07_cost}}
    \\
    \subfloat[$L_\infty$ accuracy]{\includegraphics[height=31mm]{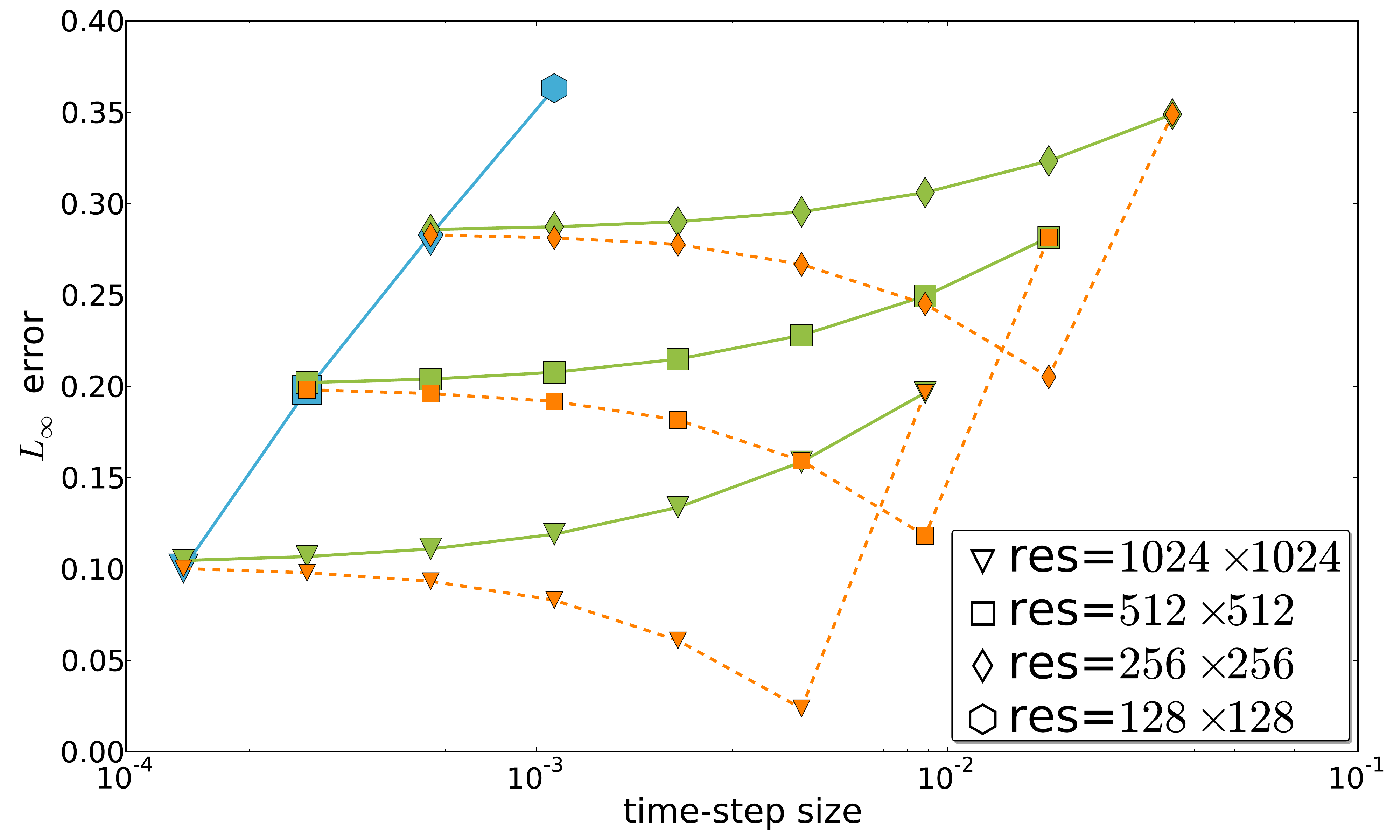}\label{fig:test07_err_inf}}
    \;
    \subfloat[$L_\infty$ cost]{\includegraphics[height=31mm]{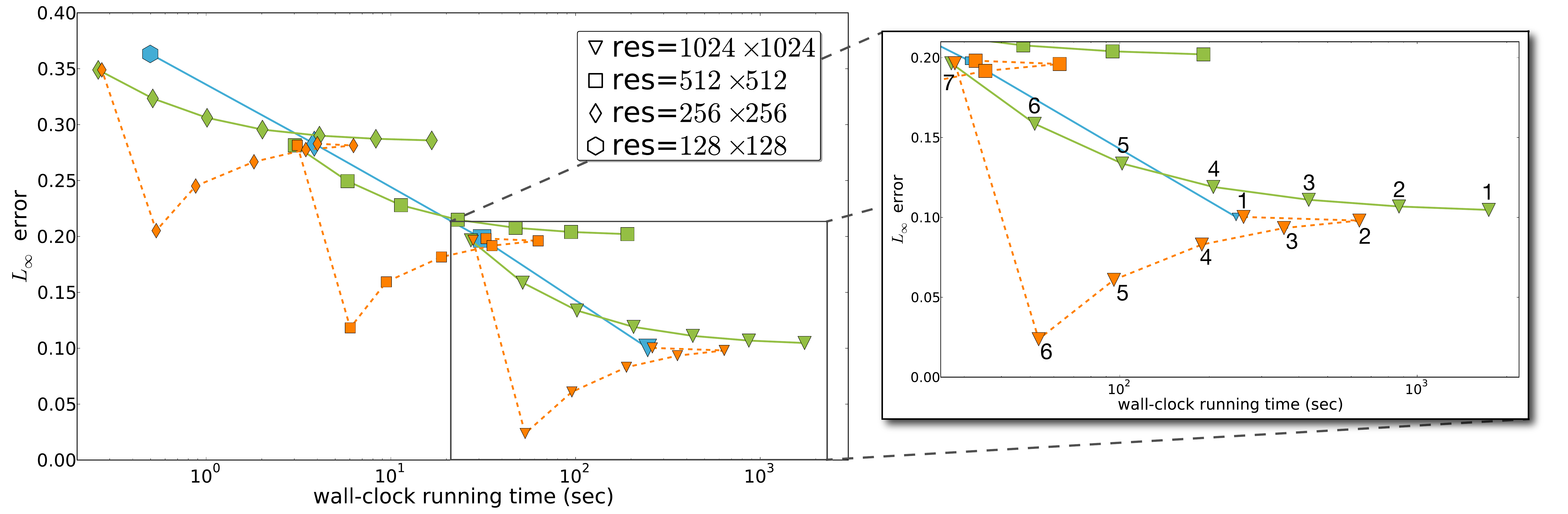}\label{fig:test07_cost_inf}}
    \caption{ {\bf Accuracy and computation cost of experiment 4}\label{fig:test07}}
\end{figure}
We consider an equation with a highly oscillatory speed function and zero boundary/terminal conditions:
\begin{equation}
\begin{split}
    u_t - \left(0.1 + 4.9\sin(\pi t)^2\sin(8\pi x)^{16}\sin(8\pi y)^{16}\right)\left|\nabla u\right| = -1,\quad & \xx\in\Omega, t\in[0,T] \\
u(\xx,t) = 0,\quad & \xx\in\partial\Omega,t\in[0,T] \textrm{ or }\xx\in\Omega, t=T
\end{split}
\end{equation}
We choose a large terminal time $T=4$ to ensure that all characteristics
passing through the time slice $t=0$ escape through $\partial\Omega$ before the terminal time $T$.
Since the analytic solution formula is not available, we first use the explicit method on a
very fine spatial grid ($2048\times2048$) to compute the ``ground truth'', which is then used to
approximate the accuracy in all experiments. The level-sets of the solution
are shown in Figure~\ref{fig:test07_ret} for three different time slices.

In this example, the speed function is periodic in $\Omega\times[0,T]$. 
It varies from $0.1$ to $5.0$, but due to the high exponents, the speed $f$ is fairly low
on most of $\Omega$ in every time slice. The large ratio between maximal and average speeds yields
very restrictive CFL conditions for the explicit methods. As a result, the implicit and hybrid methods have a clear
advantage; see Figure~\ref{fig:test07}.

The accuracy of hybrid methods here is reminiscent of the linear examples in Figure~\ref{fig:hybrid_vs_imp_vs_exp}.
When using the ``global CFL-specified'' timestep, the hybrid method defaults to the explicit update on the entire $\Omega$.
Similarly, for really large timesteps the hybrid method becomes equivalent to the implicit. For intermediate timesteps,
the hybrid's accuracy is actually better than that of both the pure (explicit and implicit) alternatives.
This happens since the CFL is locally
satisfied on most of $\Omega$, enabling explicit methods with larger timestep at those locations and thus lowering
the numerical diffusion.




Figures~\ref{fig:test07_cost} and~\ref{fig:test07_cost_inf} confirm the higher efficiency of the implicit
and hybrid methods, but also show the non-monotonic scaling of performance for the latter. The insets show
zoomed-in versions of performance curves for both methods on the $1024\times1024$ grid.
For the smallest timestep (Marker 1), the hybrid produces the same output as the explicit but with a small additional
overhead (to check if the CFL locally holds at every node).
The first timestep increase (Marker 2) leads the hybrid to use quadratic updates on a part of the domain
with an additional overhead of setting up the Fast Marching Method in each time slice. This has the net effect of
increasing the computational time, despite the fact that the number of time slices is twice smaller. As we continue to
increase the timestep, this additional overhead becomes negligible in comparison to the time slice savings, resulting
in higher efficiency (Markers 3-6). Finally, when the quadratic update is used everywhere, the hybrid's running time
becomes essentially the same as that of the implicit method (Marker 7).

We note that this example also exhibits shock lines common in viscosity solutions of HJ PDEs.
But since the characteristics run into the shock lines, their presence does not impact the performance of
implicit and hybrid methods; see also the $L_\infty$ error analysis in Figures~\ref{fig:test07_err_inf} and~\ref{fig:test07_cost_inf}.





\section{Conclusions}
\label{s:conclusions}
We have introduced new unconditionally stable methods for time-dependent
Hamilton-Jacobi-Bellman PDEs, where each next time slice is computed in
$O(M \log M)$ operations.  We have applied this approach to
fixed-horizon isotropic control problems and implemented a modified version
of Fast Marching Method to enable the time marching.
We illustrated our approach using a number of numerical examples and showed
that in many of them the new implicit and hybrid schemes outperform the explicit
time-marching method because the CFL stability conditions are found
to be more restrictive than the accuracy requirements.  This is particularly
true for speed/cost functions significantly varying on $\domain$ in each time-slice.

Even though we did not test this, we expect that other fast methods
for solving equation \eqref{eq:upwind_discr_Eikonal_from_boundary} could also be used to perform implicit
time-marching in \eqref{eq:implicit_discr_static}.  In particular, it would be interesting to
compare the efficiency of our fast marching based implementation
with that of fast sweeping based algorithms
(e.g., \cite{Zhao, TsaiChengOsherZhao, KaoOsherQian})
and other ``fast iterative'' methods
(e.g., \cite{PolyBerTsi, BorRasch, JeongWhitaker, Renzi})
in this context.

All of our examples considered in this paper were isotropic,
but the anisotropic problems could be handled similarly by using
Ordered Upwind Method \cite{SethVlad2, SethVlad3} instead of Fast Marching.
The extension to non-convex Hamiltonians (and the corresponding
fixed-horizon differential games) would require a fast solver for
the corresponding static problem.  While some such fast marching-type
algorithms are available \cite{CristianiFalcone} their computational complexity has not
been sufficiently investigated so far.  Similarly, the existing
fast sweeping-type algorithm for non-convex problems \cite{KaoOsherQian}
involves significant additional artificial viscosity, which is likely to
impact its efficiency in the context of implicit methods for time-dependent
problems.


\vspace{.2in}
\noindent
{\bf{Acknowledgments:}}
This research was supported in part by the National Science Foundation
grant DMS-1016150.
The second author's research is also supported by Columbia young faculty startup fund.

\end{document}